%% file: 0_main.tex
\documentclass[final,authoryear,5p,times,twocolumn]{elsarticle}

\usepackage{amssymb}
\usepackage{url}

\usepackage{graphicx}
\usepackage{psfrag}
\usepackage{amsmath}
\usepackage{amssymb}
\usepackage{comment}
\usepackage{color}
\usepackage{subfigure}
\usepackage{amsthm}
\usepackage[active]{srcltx}
\usepackage{threeparttable}% tables with footnotes
\usepackage{dcolumn}
\usepackage{subeqnarray}
\usepackage{subfigure}
\usepackage{subfigmat}
\usepackage{gensymb}
\usepackage{booktabs}
\usepackage{subfigure}
\usepackage{soul}

\usepackage{pgfplots}
\usepackage{tikz,flowchart}
\usetikzlibrary{shapes,arrows,chains,shapes.multipart,shapes.geometric}
\usetikzlibrary{plotmarks}
\usepgfplotslibrary{patchplots}
\usepackage{grffile}

\usepackage{natbib}

\def\eup{\,\hbox{\raise .36em\hbox to0pt{\vrule height0.5pt% 
width.55em depth0pt\hss}% 
\raise .25em\hbox to0pt{\vrule height0.5pt width.5em% 
depth0pt\hss}\hskip.02em\sf C}}

\hypersetup{
    bookmarks=true,         % show bookmarks bar?
    unicode=false,          % non-Latin characters in AcrobatÕs bookmarks
    pdftoolbar=true,        % show AcrobatÕs toolbar?
    pdfmenubar=true,        % show AcrobatÕs menu?
    pdffitwindow=false,     % window fit to page when opened
    pdfstartview={FitH},    % fits the width of the page to the window
    pdftitle={My title},    % title
    pdfauthor={Author},     % author
    pdfsubject={Subject},   % subject of the document
    pdfcreator={Creator},   % creator of the document
    pdfproducer={Producer}, % producer of the document
    pdfkeywords={keyword1} {key2} {key3}, % list of keywords
    pdfnewwindow=true,      % links in new window
    linkcolor=red,          % color of internal links
    citecolor=green,        % color of links to bibliography
    filecolor=magenta,      % color of file links
    urlcolor=cyan,           % color of external links
hyperindex=true,
pdfdisplaydoctitle=true,
hyperindex=true,
linktocpage=true,
}
\newcommand{\R}{\mathbb{R}}

\newcommand{\m}{\mathcal}
\usepackage{subfigure}
\usepackage{subfigmat}
\usepackage{mathtools}
\usepackage{algorithmicx}
\usepackage{algorithm}
\usepackage{algpseudocode} % for pseudocode
\usepackage{amssymb, latexsym, mathrsfs}
\usepackage{multicol}

\algnewcommand\algorithmicinput{\textbf{Input:}}
\algnewcommand\INPUT{\item[\algorithmicinput]}

\algnewcommand\algorithmicoutput{\textbf{Output:}}
\algnewcommand\OUTPUT{\item[\algorithmicoutput]}

\usepackage{amsthm}

\theoremstyle{definition}

\newtheorem*{definition*}{Definition}

\usepackage{float}
\floatstyle{plaintop}
\restylefloat{table}

\usepackage{nomencl}

\journal{\scriptsize{Elsevier Aerospace Science and Technology}}
\begin{document}

\begin{frontmatter}

%\title{Maximizing Safety of Aircraft Trajectories in Stochastic Weather}
%\title{Aircraft Trajectory Planning to Maximize Safety in the Presence of Thunderstorms}
\title{On Maximizing Safety in Stochastic Aircraft Trajectory Planning \\ with Uncertain Thunderstorm Development}
%\title{On Maximizing Safety in Stochastic Aircraft Trajectory Planning under the Presence of Thunderstorms}

\author[Hentzen,Kamgarpour,Soler,Gonzalez]{Daniel Hentzen, Maryam Kamgarpour, Manuel Soler, Daniel Gonz\'alez-Arribas}
\address[Hentzen]{Automatic Control Lab, ETH Zurich, Zurich, Switzerland. e-mail: dhentzen@ethz.ch%
}% <-this % stops a space
\address[Kamgarpour]{Automatic Control Lab, ETH Zurich, Zurich, Switzerland. e-mail: mkamgar@control.ee.ethz.ch}
\address[Soler]{Department of Bioengineering and Aerospace Engineering, Universidad Carlos III de Madrid, Legan\'es, Spain. e-mail: masolera@ing.uc3m.es}
\address[Gonzalez]{Department of Bioengineering and Aerospace Engineering, Universidad Carlos III de Madrid, Legan\'es, Spain. e-mail: daniel.gonzalez.arribas@uc3m.es}

\begin{abstract}
Dealing with meteorological uncertainty poses a major challenge in air traffic management (ATM). Convective weather (commonly referred to as storms or thunderstorms) in particular represents a significant safety hazard that is responsible for one quarter of weather-related ATM delays in the US. With commercial air traffic on the rise and the risk of potentially critical capacity bottlenecks looming, it is vital that future trajectory planning tools are able to account for meteorological uncertainty. We propose an approach to model the uncertainty inherent to forecasts of convective weather regions using statistical analysis of state-of-the-art forecast data. The developed stochastic storm model is tailored for use in an optimal control algorithm that maximizes the probability of reaching a waypoint while avoiding hazardous storm regions. Both the aircraft and the thunderstorms are modeled stochastically.
The performance of the approach is illustrated and validated through simulated case studies based on recent nowcast data and storm observations.
\end{abstract}


\end{frontmatter}

% \linenumbers

\input{1_intro.tex}

\input{2_stoch_model.tex}
\input{3_stoch_reach.tex}
\input{4_case_studies.tex}
\input{5_conclusion.tex}

%%%%%%%%%%%%%%%%%%%%%%%%%%%%%%%%%%%%%%%%%%%%%%%%
\section*{Acknowledgments}
%%%%%%%%%%%%%%%%%%%%%%%%%%%%%%%%%%%%%%%%%%%%%%%%
We would like to thank the \textit{Agencia Estatal de Meteorolog\'ia} (AEMET) for kindly providing the nowcast data employed in this research. M. Kamgarpour thanks Swiss National Fund 200021\_172782 for their partial support.




%\appendix
%
%%\input{APPA_reformulation}
%%\input{APPB_5thdegree}
%
%\newpage
%%\input{APPC_MINLPapproach}
%
%\clearpage

%% References


%\bibliographystyle{elsarticle-num}
%\bibliographystyle{apa-good}
\bibliographystyle{model2-names}
\bibliography{library.bib}

%% not want to use model2-names.bst.



%% \bibitem must have one of the following forms:
%%   \bibitem[Jones et al.(1990)]{key}...
%%   \bibitem[Jones et al.(1990)Jones, Baker, and Williams]{key}...
%%   \bibitem[Jones et al., 1990]{key}...
%%   \bibitem[\protect\citeauthoryear{Jones, Baker, and Williams}{Jones
%%       et al.}{1990}]{key}...
%%   \bibitem[\protect\citeauthoryear{Jones et al.}{1990}]{key}...
%%   \bibitem[\protect\astroncite{Jones et al.}{1990}]{key}...
%%   \bibitem[\protect\citename{Jones et al., }1990]{key}...
%%   \harvarditem[Jones et al.]{Jones, Baker, and Williams}{1990}{key}...
%%


% \end{thebibliography}

\end{document}

%% file: 1_intro.tex
% !TeX spellcheck = en_US
% !TeX encoding = UTF-8
% !TeX root = 0_main.tex

%%%%%%%%%%%%%%%%%%%%%%%%%%%%%%%%%%%%%%%%%%%%%%%%
\section{Introduction}
\label{sect:introduction}
%%%%%%%%%%%%%%%%%%%%%%%%%%%%%%%%%%%%%%%%%%%%%%%%

% motivation: weather delays/disturbances in aviation

A challenging aspect of tactical aircraft trajectory planning\footnote{By \textit{tactical} we refer to 10 minutes to 1 hour prior to the potential hazard encounter.} is the avoidance of hazardous convective weather regions, more commonly referred to as storms or thunderstorms\footnote{Throughout this paper we will use the more general term \textit{thunderstorm} to describe convective weather regions.}. The inherently uncertain nature of thunderstorms causes major safety risks. Strong conflicting up- and downdrafts lead to heavy turbulence. Hail, severe icing and lightning can also inflict significant damage to aircraft equipment and windshields. In addition to increased safety risks, thunderstorms are also a leading cause of reduced time (\cite{NationalResearchCouncil}) and cost (\cite{Michalek2009a}) efficiencies. From 2008 to 2013, inclement weather caused 69\% of system-impacting delays (delays greater than  15 minutes), as recorded in the \textit{OPSNET Standard ``Delay by Cause" Reports}\footnote{{NextGen Weather - Weather Delays} \url{https://www.faa.gov/nextgen/programs/weather/faq/}}. Within those weather delays, thunderstorms emerging from atmospheric instabilities were responsible for around 25\%, turning them into the leading cause of flight delays in the US airspace. 

% existing approaches to deal with convective weather
Today, commonly used approaches to plan trajectories in the presence of thunderstorms are suboptimal and prone to human error. Thunderstorms are typically detected by satellite images or weather radar signals, which are processed into graphical images. The identified \textit{thunderstorm cells} are then classified according to severity. Regions of airspace that are considered to pose a significant safety hazard are recommended as unsafe-to-fly regions to air traffic controllers (ATCO) and pilots (\cite{Krozel2007}). 
 
In many cases these regions are avoided by manual trajectory modifications. This involves increased workload and potential loss of situational awareness for both ATCOs and pilots, jeopardizing safety and air traffic capacity.

% recent efforts
One of the key concepts of future Air Traffic Management (ATM) systems is the so-called \textit{trajectory-based operations} (TBO) described in the SESAR (\cite{sesar}) or NextGen (\cite{nextgen}) ATM master plans. The TBO concept will give aircraft the autonomy to update their trajectories to satisfy business needs, such as fuel or flight time savings, or to avoid potential conflicts or hazardous weather regions. Robust planning that accounts for the uncertainty associated to aircraft motion and thunderstorm development is paramount for a successful transition to TBO.

  In \cite{mcnally2015dynamic}, the authors reported two years of operational testing of a tool coined {\it Dynamic Weather Routes} at American Airlines. The tool is a ground-based real-time trajectory automation system that continuously analyzes en-route aircraft. It computes simple corrections to flight plans to avoid thunderstorms while ensuring conflict avoidance, respecting airspace constraints and minimizing delay. However, the tool neither accounts for weather uncertainty, nor does it use advanced control and optimization algorithms, thus leaving room for improvement in the system's efficiency, safety, and capacity.

% p3 the need for accurate forecast to design automation
The lack of probabilistic forecast data is one of the main bottlenecks in the development and full-scale deployment of automated trajectory planning tools that account for thunderstorm uncertainty. Today's state-of-the-art meteorological forecast products include current observations of the thunderstorm cells as well as predictions for their locations at discrete forecast horizons. A commonly used approach relies on deterministic short term forecasts, so-called nowcasts. They provide information on thunderstorm cell positions for up to a few hours ahead (\cite{wilson1998nowcasting}). Examples for such nowcast systems include the satellite data-based cloud tracker Cb-TRAM (Cumulonimbus Tracking and Monitoring) \cite{zinner2008cb}; and the radar-based Rad-TRAM (Radar Tracking and Monitoring), developed by \cite{kober2009tracking}.
Nevertheless, these forecast products do not provide any uncertainty quantification.

The integration of stochastic information in meteorological forecast products has indeed proven difficult, mainly due to the short lifespan of convective phenomena (around 30 minutes) and an imperfect knowledge of the atmosphere at macroscopic scale. Recent efforts in this direction include, e.g., \cite{mecikalski2015probabilistic, osinski2017short}. The former presented two statistical learning approaches that incorporate numerical weather prediction (NWP) input within satellite images to produce probabilistic forecasts. \cite{osinski2017short} presented an experimental product developed by Meteo-France. It combines the numerical weather prediction model AROME-NWC with simulated radar images and post-processing. An ensemble of time-lagged forecasts with thunderstorm development information is then produced. Besides these efforts in combining data-driven methods and NWP, pure data-driven methods can be also used to model uncertainties in forecast products through comparison of actual instances of thunderstorms and their predictions. 

The authors in \cite{Radhakrishna2012} investigated statistical properties of thunderstorm cells' growth and decay, as well as their scale dependence and predictability. They concluded that forecast accuracy of deterministic, radar-based nowcasts can be improved by the extrapolation of growth and decay for a few hours. 
In \cite{Sauer2014}, the authors analyzed the error of radar-based nowcasts and their increasing uncertainty with time by determining spatial deviations of the nowcast predictions compared to the observed thunderstorm realizations. The authors argued that thunderstorm cells could be enlarged by an uncertainty margin that is dependent on the desired probability of avoiding the respective thunderstorm cell. 
In \cite{zhang2013probabilistic,zhang2017probabilistic}, the authors introduced a probabilistic analysis of weather forecasting. The analysis is based on ensemble weather forecasting data provided by the Center for Analysis and Prediction of Storm (CAPS) \cite{jung2008assimilation}. %These findings serve as a basis for the data-driven methodology introduced in this paper for  modelling  storms' inherent stochasticity.
Additionally, \cite{zhang2017probabilistic} presented a path planning algorithm to avoid the stochastic hazard, including a mission risk analysis. Nevertheless, additional efforts are needed in stochastic thunderstorm modeling and aircraft trajectory planning, e.g., considering the additional uncertainty introduced by the aircraft dynamic system.

% p6 past work on addressing trajectory design in uncertain forecast/obstacles
Several methods to address aircraft path planning subject to uncertainties in system dynamics and obstacles (such as unsafe-to-fly regions) have been proposed. The problem of routing an aircraft to a target while avoiding uncertain hazardous weather zones in the presence of wind disturbance can be cast as the optimal control of a stochastic dynamic system through probabilistic obstacles, such as unsafe-to-fly regions. The following challenges arise: First, both the aircraft dynamics and the thunderstorm cells are nonlinear and non-convex. Second, the use of stochastic models is required, resulting in a theoretically and computationally challenging \textit{stochastic optimal control} problem. 

A possible approach is to bound the uncertain obstacles with deterministic sets. 
The authors in \cite{Kuwata2005} presented a receding horizon control strategy for aircraft trajectory planning in environments with unknown but bounded disturbances. The authors in \cite{Kamgarpour2010a} proposed a receding horizon framework for designing safe aircraft trajectories. They took into account a dynamic weather forecast product to determine deterministic bounds for the thunderstorm regions to be avoided. While this set-bounded approach offers a high degree of safety, the solutions often degrade overall performance due to high conservatism arising from the protection against low probability uncertainty outliers. More importantly, the optimization problem can quickly become infeasible, for example in a scenario featuring a high number of  obstacles. Such scenarios call for a solution that minimizes the risk of failure, instead of trying to eliminate it entirely.
In \cite{Liu2014}, a stochastic optimal control framework for mid-air conflict resolution was presented that could incorporate uncertainties in both aircraft and wind dynamics. However, in the studied aircraft-weather conflict, the weather hazard was considered deterministic. In \cite{Nilim2002}, the authors addressed the problem of generating aircraft trajectories using a Markov decision problem where the evolution of the thunderstorms was modeled as a Markov chain. The transition probabilities were extracted from historical data using maximum-likelihood estimators. However, these probabilities were assumed stationary. Moreover, delay was used as the objective to be minimized. Given the uncertainties in thunderstorm development, an important objective is to minimize the probability of intersecting a hazardous weather region. 

A suitable approach to generate trajectories that reach a target set while minimizing the probability of intersecting probabilistic unsafe-to-fly regions is to formulate a \textit{stochastic reach-avoid problem}. The authors in \cite{Summers2011} considered a stochastic reach-avoid problem with a time-varying stochastic obstacle set. The  problem was reformulated as a finite horizon stochastic optimal control problem and solved with Bellman's Dynamic Programming principle. We proposed this method as a feasible approach to aircraft trajectory planning in the presence of a single hypothetical thunderstorm, modeled by a stochastic ellipse \cite{Hentzen2017}.  However, in this preliminary study, we assumed a priori knowledge of the stochastic unsafe-to-fly regions. No forecast data was incorporated to define stochastic unsafe-to-fly regions. Furthermore, the approach did not address trajectory planning in the presence of multiple thunderstorm regions.  

%p7: goals
The contribution of the paper is twofold: First, we propose a method for modeling thunderstorm in a stochastic manner. The method accounts for uncertainty in both the movement and, crucially, in the growth or decay of the thunderstorm cells. It is based on statistical analysis of  historical nowcast data and its comparison with actual thunderstorm observations. Second, we combine the developed stochastic thunderstorm model with an optimal trajectory planning algorithm based on the stochastic reach-avoid methodology by \cite{Summers2011}. The resulting aircraft trajectories  maximize the probability of reaching a given waypoint while avoiding multiple thunderstorms, taking into account uncertainties in both the system dynamics (due to wind disturbance) and the unsafe-to-fly regions (due to uncertainties in nowcast data on thunderstorms). The resulting trajectories are validated against actual realizations of thunderstorms. 

The paper is structured as follows: The proposed stochastic thunderstorm model derived from nowcast data is presented in Section~\ref{sect:stoch_mod}. The stochastic reach-avoid framework, in which the aircraft trajectory generation is cast as a stochastic optimal control problem, is presented in Section~\ref{sect:stoch_reach}. The trajectory planner, resulting from the integration of the stochastic thunderstorm model and the optimal control algorithm into one framework, is applied to two flight scenarios using historical forecast data in Section~\ref{sect:case_studies}. Conclusions and future research directions are discussed in Section~\ref{sect:conclusion}.

%% file: 2_stoch_model.tex
% !TeX spellcheck = en_US
% !TeX encoding = UTF-8
% !TeX root = 0_main.tex

%%%%%%%%%%%%%%%%%%%%%%%%%%%%%%%%%%%%%%%%%%%%%%%%
\section{Stochastic Thunderstorm Modeling}
\label{sect:stoch_mod}
%%%%%%%%%%%%%%%%%%%%%%%%%%%%%%%%%%%%%%%%%%%%%%%%

\begin{figure}
	\centering
	\includegraphics[width=0.8 \linewidth]{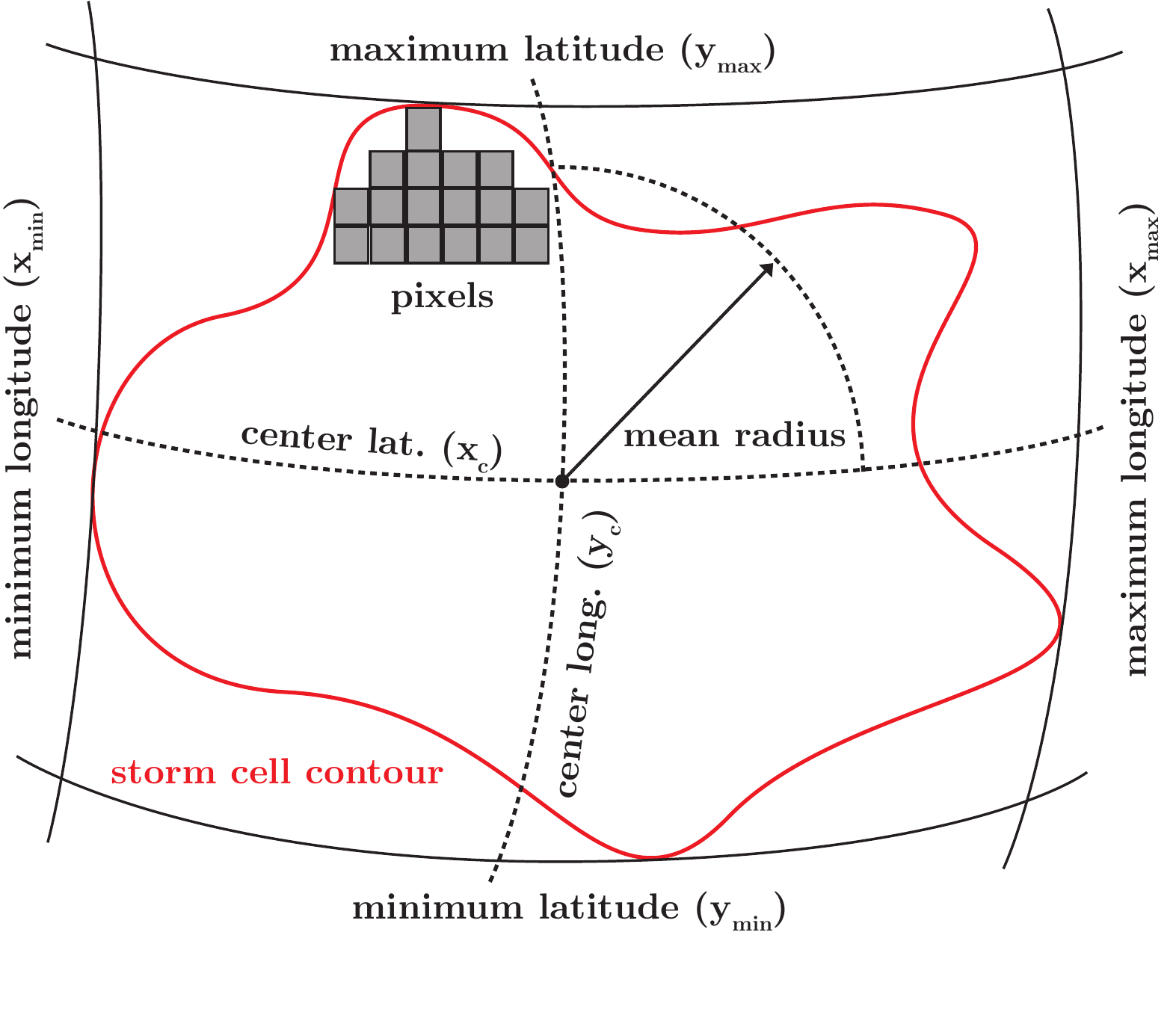}
	\caption{Thunderstorm cell representation}
	\label{fig:storm_cell}
\end{figure}

The unsafe-to-fly regions to be avoided by the aircraft are hazardous convective weather regions, referred to as \textit{thunderstorm cells}. The actual position and size of these regions, as well as an extrapolation of their position at different forecast horizons, are provided in the form of \textit{nowcasts}. A nowcast file contains observations of the current state of the hazardous thunderstorm cells that are present in a particular portion of the airspace (here denoted by the 2D flight level $X \times Y$). These observations include, among others, the position of each thunderstorm cell center (in longitude, latitude coordinates), the extremities of each thunderstorm cell in all four directions (given in longitude, latitude), the size of each thunderstorm cell (expressed in number of ``pixels''), a mean radius, as well as the speed and direction of each thunderstorm cell. Figure \ref{fig:storm_cell} depicts  these thunderstorm cell characteristics. A nowcast file also provides a series of extrapolated predictions at discrete forecast horizons. However, these predictions are limited to the position of the thunderstorm cell centers. No other thunderstorm cell parameters are predicted. Every 10 minutes, a new nowcast file is published. Each nowcast file gives predictions for 6 forecast horizons from 10 to 60 minutes with a 10 minute spacing. Compared to the eventual observations of the thunderstorms, these forecasts feature significant prediction errors, mainly due to the nonlinear growth and decay of rainfall field size, which according to \cite{Browning1982} account for up to a quarter of the forecast error. This motivates the use of a stochastic model to represent the thunderstorm obstacles.

\begin{figure}
	\centering
	\includegraphics[width=0.8\linewidth]{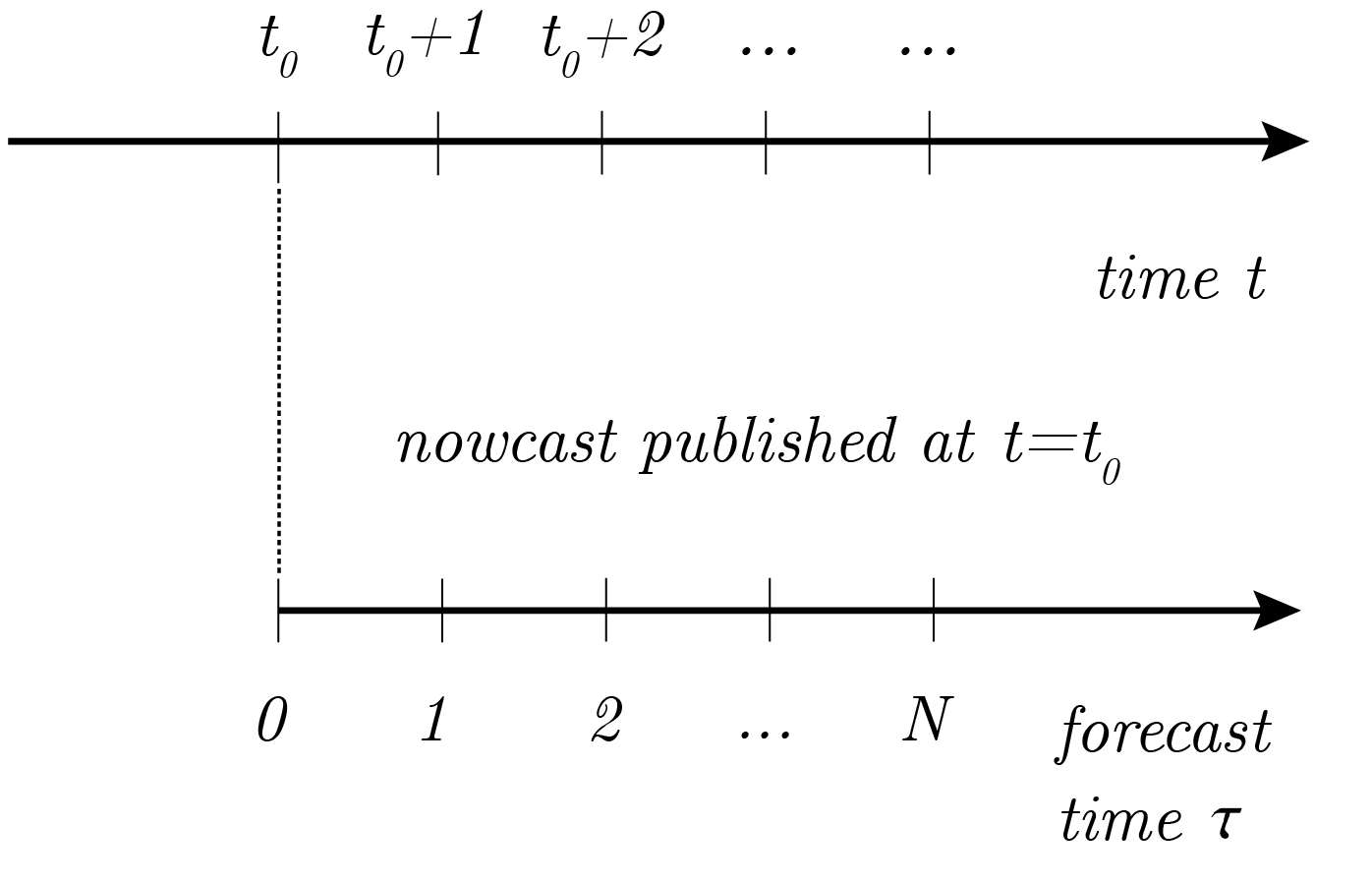}
	\caption{Nowcast structure.}
	\label{fig:timeline}
\end{figure}

We consider a constant altitude flight, and thus, a two-dimensional domain. We aim to evaluate the probability that any discrete region of the flight level will be affected by hazardous thunderstorm at different forecast horizons $\tau=1, 2\dots,N$, corresponding to the forecast times 10, 20, $\dots$, 60 minutes and interpolated time instances in between. Figure \ref{fig:timeline} illustrates the nowcast structure and notation. Let $\Xi_{t_0+\tau|t_0}$ be the set of 2D \textst{obstacles} regions capturing the thunderstorm forecast for time $t_0+\tau$ in the nowcast published at time $t_0$. For the sake of brevity, we define $\square_{\tau} \coloneqq \square_{t_0+\tau|t_0}$ throughout the paper. Let $X \times Y$ denote the 2D airspace under consideration. At each time instance $t$, we assign to each discretized coordinate $(x,y) \in X \times Y$ a thunderstorm probability, that is, the probability that the position $(x,y)$ is in a thunderstorm cell. For this purpose we use the function $p_{\Xi_{\tau}}(x,y): X\times Y \mapsto [0,1]$. To compute $p_{\Xi_{\tau}}(x,y)$ for $\tau=1,\dots,N$, starting from a known and deterministic initial obstacle set $\Xi_{0}$, a stochastic model that describes the evolution of the thunderstorm is needed. 

All information on thunderstorm position, size and movement is drawn from the available nowcasts. The analysis herein relies on a data sample of 720 nowcast files collected between December 15th and December 19th 2016, provided by the Spanish meteorological agency, \textit{Agencia Estatal de Meteorolog\'ia} (AEMET). AEMET's YRADAR application produces both a 2D and a 3D analysis. Due to the simplifying assumption of constant altitude, the following analysis will be constrained to the 2D dataset. Due to the short average lifespan of thunderstorm cells, a significant proportion of detected cells disappear after the 40 minute horizon, rendering the data for the 50 and 60 minute horizons sparse. We will therefore restrict the analysis to a forecast horizon up to 40 minutes. The most relevant nowcast parameters used for the stochastic characterization of the storm cells are presented in Table \ref{tab:parameters_nowcast}.
\begin{table}[b!]
\centering
\begin{scriptsize}
\begin{tabular}{p{2.95cm} p{4.5cm}}
	\hline
	\textbf{Nowcast parameter} & \textbf{Description}\\
	\hline
	NUM & ID of the cell\\
	NUPIX & Number of pixels of the cell\\
	LATCEN, LONCEN & Latitude/longitude of cell's center\\
	RADIOE & Effective radius of the cell\\
	LATNOR, LATSUR & Limits of the cell in latitude\\
	LONOES, LONEST & Limits of the cell in longitude\\
	LAT10, LON10 & Forecast for center after 10 min.\\
	LAT20, LON20 & Forecast for center after 20 min.\\
	\vdots & \vdots\\
	LAT60, LON60 & Forecast for center after 60 min.\\
	DIRN & Direction of displacement of the storm cell\\
	V(KMH) & Velocity of displacement of the storm cell\\
	\hline
\end{tabular}
\end{scriptsize}
\caption{Overview of relevant nowcast parameters.}
\label{tab:parameters_nowcast} 
\end{table}

We represent the thunderstorm cells using two states:
\begin{enumerate}
	\item the position of the cell center $(x_t^c,y_t^c)$ at time $t$,
	\item the position of the thunderstorm cell's lateral and longitudinal extremities $(x^{(E)}_t, x^{(W)}_t,y^{(N)}_t,y^{(S)}_t)$ at time $t$ .
\end{enumerate}
 While no information on the thunderstorm cell's shape is available, the four extremities enable us to formulate a rectangular bound of the thunderstorm cell. At time $t=t_0$ ($t_0$ being the time at which the nowcast is published), we have deterministic knowledge (assumed exact) of all six state variables. For times $t=t_0 + \tau$ $(\tau =1,\dots,N)$, the nowcast files provide a prediction of the thunderstorm cell centers but no information on the thunderstorm cell extremities. To account for the prediction error on the center position (increasingly significant with growing forecast horizon) and the lack of extremity predictions, we aim to formulate a stochastic dynamic model, which will allow us to predict the evolution of the thunderstorm cell state in a probabilistic fashion. By using this model in combination with a sampling-based method, we will be able to compute the probability of thunderstorm at each $(x_t,y_t) \in X\times Y$ for each time instance $t = t_0 + \tau$, that is, $p_{\Xi_\tau}(x,y)$.

\subsection{Modeling Stochastic Thunderstorm Cell Dynamics}
\label{subsec:StormDyn}

Let $\mathbf{c}_t = (x^c_t,y^c_t) \in \R^2$ denote the observed center of the thunderstorm cell at time $t$, where $c$ denotes the center. Let $\mathbf{c}_{\tau}$ denote the forecast made at time $t_0$ for the thunderstorm cell center $\mathbf{c}$ at time $t_0+\tau$. Recall that $\tau=1,2,\dots,N$ is the forecast horizon corresponding to the forecast times 10, 20, $\dots$, 40 minutes and interpolations in between.

% Should we limit the forecast horizon to i_max instead of t?

The forecast dynamics for the (observed) thunderstorm cell center position is given by $\mathbf{c}_{t_0+\tau} = \mathbf{c}_{\tau} + \eta_{\tau}$, 
where $\eta_\tau = (\eta_{\tau,x}, \eta_{\tau,y})$ denotes the prediction error associated with the forecast horizon $\tau$ in x and y directions. The prediction error grows with the length of the forecast horizon (\cite{Sauer2014,Radhakrishna2012}). It may also be dependent on current meteorological conditions, such as strong wind in a particular direction. In the developed model, this error is accounted for by sampling it from a parameterized probability distribution and adding it to the nowcast prediction. The distribution type is determined by analyzing historical predictions and comparing them with the observations. Once an appropriate distribution type has been determined, the parameters of the distribution are determined and updated in real time using maximum likelihood estimation.

\begin{figure*}[t!]
\includegraphics[width=1\linewidth]{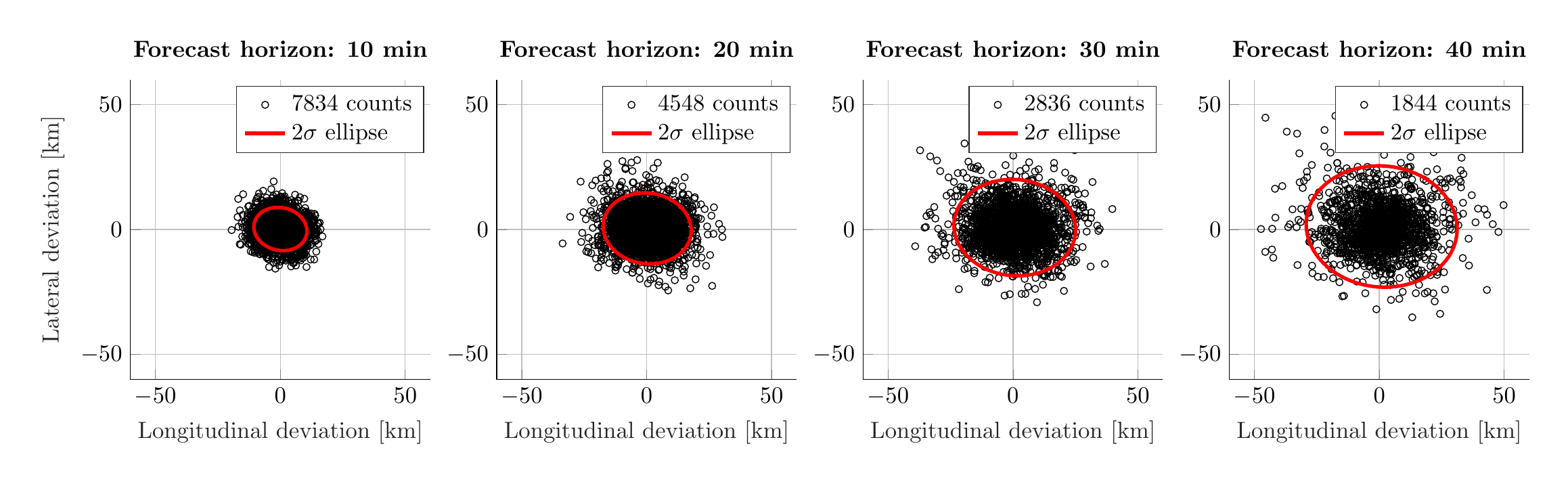}
\vspace*{-0.5cm}
\caption{Longitudinal and lateral deviation of forecast thunderstorm cell center positions to observed realizations over growing forecast horizon.}
\label{fig:scatterhist_center}
\end{figure*}

The forecast error at forecast horizon $\tau$ is given by $\eta_{\tau} = \mathbf{c}_{t_0+\tau} - \mathbf{c}_{\tau}$. Figure \ref{fig:scatterhist_center} shows combined scatter plots and histograms illustrating the deviation in longitude and latitude of the forecast center positions compared to the observed realizations in the training set for forecast horizons 10, 20, 30 and 40 minutes. Using Matlab's \texttt{allfitdist} function, all valid parametric probability distributions are fitted over the empirical data and compared for goodness-of-fit, as quantified by the Bayesian information criterion (BIC), a criterion for selecting a model among a finite set of models. The \textit{logistic} distribution offers the best fit for the available data (see \cite[Chapter 3]{HentzenThesis} for a detailed discussion, where logistic, tlocationscale, normal and generalized extreme value distribution fittings are presented, compared, and discussed). In particular, it is able to better account for the heavier tails compared to a normal distribution.

The probability density function (PDF) $f_{\mathcal{L}}$ of the \textit{logistic} distribution $\m{L}$ is given by
\begin{equation*}
f_{\mathcal{L}}(\eta;m,s) = \frac{1}{4s} \text{sech}^2\bigg(\frac{\eta-m}{2s}\bigg) = \frac{e^{-\frac{\eta-m}{s}}}{s\big(1+e^{-\frac{\eta-m}{s}}\big)^2},
\end{equation*}
where $\eta \in \mathbb{R}$ is the forecast error, $m$ is the mean of the distribution and $s\geq0$ is the \textit{scale} or \textit{dispersion} parameter. The standard deviation $\sigma$ is expressed by $\sigma = \frac{s \pi}{\sqrt{3}}$.  The distribution's parameters are determined using maximum likelihood estimation with a defined number of historical data points. 
Based on the studied data samples, we conclude:
\begin{enumerate}
	\item The thunderstorm cell center position prediction error in x and y directions $\eta_{\tau,x}$ and $\eta_{\tau,y}$ at forecast horizon $\tau$ can be approximated by sampling from the logistic PDFs $\m{L}_{x}^c$ and $\m{L}_{y}^c$ associated to the forecast horizon $\tau$: $\eta_{\tau,x} \sim \m{L}_{x}^c(\tau)$, $\eta_{\tau,y} \sim \m{L}_{y}^c(\tau)$.
	\item The parameters of these distributions (mean and scale) at different forecast horizons are computed using maximum likelihood estimation. 
\end{enumerate}

 As the available nowcasts offer no forecast on thunderstorm cell size, past observations have to be analyzed to identify and extrapolate a growth or decay trend.

In the following analysis, growth and decay are expressed in terms of the change in lateral and longitudinal size of a single thunderstorm cell. Longitudinal size (width $w$) and lateral size (height $h$) at time $t$ are defined using the nowcast extremity data:
\begin{align*}
\begin{split}
w_t &= x^{(E)}_t - x^{(W)}_t, \quad h_t = y^{(N)}_t - y^{(S)}_t.
\end{split}
\end{align*}
The dynamics of the thunderstorm cell size are given by
\begin{align*}
\begin{split}
w_{t+1} &= w_{t} + \Delta w, \quad h_{t+1} = h_{t} + \Delta h,\\
\end{split}
\end{align*}
where $\Delta w$ and $\Delta h$ denote the change in width and height, respectively. They are the stochastic parameters to be extrapolated from historical data.

Similar to the analysis of the thunderstorm cell center prediction error above, we determine an appropriate parameterized probability distribution to model $\Delta w$ and $\Delta h$. The BIC indicates that the logistic probability distribution offers a superior fit to the normal distribution (see \cite{HentzenThesis} for a detailed discussion). We observe that thunderstorm cells featuring a large surface area (expressed using the number of pixels on the radar images) grow or decay faster than smaller cells. We thus choose the parameters of the \textit{logistic} distribution as a function of the thunderstorm cell size in pixels ($s_{pix}$):
\begin{align*}
&\Delta w \sim \m{L}^{w}(s_{pix}) \quad \Delta h \sim \m{L}^{h}(s_{pix}).
\end{align*}
The scale parameter of the $\Delta h$ and $\Delta w$ distributions is assumed to grow logarithmically with thunderstorm cell size (see \cite{HentzenThesis}).

\label{para:stoch_model}
The derived stochastic model of the thunderstorm cell dynamics is given by two steps as follows.
\begin{enumerate}
	\item Sample step:
	\begin{align*}
	&\eta_{\tau,x} \sim \m{L}_{x}^c(\tau), \quad \eta_{\tau,y} \sim \m{L}_{y}^c(\tau)\\
	&\Delta w \sim \m{L}^{w}(s_{pix}), \quad \Delta h \sim \m{L}^{h}(s_{pix})
	\end{align*}
	\item Update step:
	\begin{equation*} \label{eq:cell_dyn}
	\begin{cases}
	\underbrace{\begin{bmatrix}
		x^c_{t_0+\tau}\\
		y^c_{t_0+\tau}\\
		\end{bmatrix}}_{\substack{\text{corrected centers} \\ \text{at forecast horizon $\tau$}}} = 
	\underbrace{\begin{bmatrix}
		x^{c}_{\tau}\\
		y^{c}_{\tau}\\
		\end{bmatrix}}_{\substack{\text{uncorrected centers} \\ \text{at forecast horizon $\tau$}}} +
	\underbrace{\begin{bmatrix}
		\eta_{\tau,x}\\
		\eta_{\tau,y}\\
		\end{bmatrix}}_\text{sampled error}\\
	\qquad\\
	\underbrace{\begin{bmatrix}
		x_{t_0+\tau}^{(W)}\\
		x_{t_0+\tau}^{(E)}\\
		y_{t_0+\tau}^{(S)}\\
		y_{t_0+\tau}^{(N)}\\
		\end{bmatrix}}_{\substack{\text{corrected extremities} \\ \text{at forecast horizon $\tau$}}} = 
	\begin{bmatrix}
		x^c_{t_0+\tau}\\
		x^c_{t_0+\tau}\\
		y^c_{t_0+\tau}\\
		y^c_{t_0+\tau}\\
		\end{bmatrix} +
	\underbrace{\begin{bmatrix}
		-x^c_{t_0} + x_{t_0}^{(W)}\\
		x_{t_0}^{(E)} -x^c_{t_0}\\
		-y^c_{t_0} + y_{t_0}^{(S)}\\
		y_{t_0}^{(N)} -y^c_{t_0}\\
		\end{bmatrix}}_{\substack{\text{width/height} \\ \text{at $t_0$}}}+\frac{1}{2}
	\underbrace{\begin{bmatrix}
		- \sum_{i=1}^{\tau} \Delta w(s_{pix})\\
		\sum_{i=1}^{\tau}\Delta w(s_{pix})\\
		-\sum_{i=1}^{\tau} \Delta h(s_{pix})\\
		\sum_{i=1}^{\tau} \Delta h(s_{pix})
		\end{bmatrix}}_{\substack{\text{sampled change} \\ \text{in width/height} \\ \text{over forecast horizon $\tau$}}}
	\end{cases}
	\end{equation*}
\end{enumerate}
.

\subsection{Obstacle Clustering}
\label{subsec:StochNoFly}
The dynamic thunderstorm cell model provides a sequence of predictions of thunderstorm cell positions and sizes for $t=t_0 + \tau$ with $\tau = 1, \dots, N$. However, the individual thunderstorm cells alone are not suitable obstacles  for scenarios with a high density of convective weather. The thunderstorm cells are often grouped closely together in the airspace affected by convective weather. If every thunderstorm cell were to be considered as an individual obstacle, the resulting safety-optimal trajectory would likely contain a high number of turn advisories, as individual obstacles would have to be avoided every few minutes. Although such a trajectory might be feasible in theory, it is not an option in a realistic ATM setting. Furthermore, individual thunderstorm cells are likely to split into sub-cells or merge with other neighboring or overlapping cells. It is also likely that new thunderstorm cells appear in regions with high levels of convective weather. Nowcasts  do not contain information on these processes. Considering the thunderstorm cells as individual obstacles would thus result in potentially unsafe trajectories. For the these reasons, we propose a more robust approach in which the individual thunderstorm cells predicted for forecast horizons $\tau=1,\dots, N$ by the dynamic thunderstorm cell model at time $t_0$, are clustered into larger \textst{obstacles} unsafe-to-fly regions. This clustering is updated at each forecast horizon $\tau$ to account for the movement of the individual thunderstorm cells. The clustering results in a sequence of \textst{obstacles} unsafe-to-fly regions that represent areas of the airspace that feature a high density of thunderstorm cells and a high probability of thunderstorm occurrence. We use K-means clustering for this purpose. 

At time $t_0 + \tau$, we forecast the individual thunderstorms $(\mathbf{o}^1_\tau, \mathbf{o}^2_\tau, \dots, \mathbf{o}^n_\tau)$, where each element is a $d$-dimensional real vector. The $n$ forecast thunderstorms are partitioned into $K \leq n$ clusters $\{C^1_\tau, C^2_\tau, \dots, C^K_\tau\}$ with centroids (that is, mean of points within the cluster) $\{\bar{\mathbf{o}}^1_\tau, \bar{\mathbf{o}}^2_\tau,\dots, \bar{\mathbf{o}}^K_\tau\}$ such that the variance of the distance of the observations within a cluster to the cluster's centroid is minimized:
\begin{equation*}
\arg \min_{C} \sum_{k=1}^K \sum_{\mathbf{o^i_\tau} \in C^k_\tau } \lVert \mathbf{o}^i_\tau-\bar{\mathbf{o}}^k_\tau \rVert^2_2.
\end{equation*}
Above,  $\lVert \cdot \rVert_2$ is the Euclidean $2$-norm. For our application, we consider $\mathbf{o}_\tau =(x^c_{\tau}, y^c_{\tau}, \xi_{\tau})$, where $x^c$ and $y^c$ are the $x$ and $y$ coordinates of the thunderstorm center and $\xi$ is the heading of the thunderstorm, as provided in the nowcast. This choice allows to group the thunderstorms based on their distance to each other and based on the directions in which they are heading.

Choosing the number of clusters $K$ is a hard algorithmic problem. A possible metric for the quality of the resulting clusters is the sum of squared errors (SSE), that is, the sum of the squared distance between each member of the cluster and its centroid. The SSE decreases with growing $K$, with $\text{SSE} = 0$ for $K=n$. Applying the so-called \textit{elbow method}, we choose $K$ at the \textit{elbow point}, that is, where the SSE starts decreasing at a significantly slower rate. The elbow method is heuristic and is not guaranteed to give satisfactory results in all cases. We discuss an alternative, empirical method to determine $K$ based on a parameter study with historical data in \cite{HentzenThesis}. %Alternatively, we can classify the nowcast files according to the contained number of storm cells. For each class, we compare the performance of different $K$.we use a look-up table mapping the number of storm cells present in the nowcast to the number of clusters $K$. This map is generated in a parameter study using historical data such that it maximizes forecast skill over all considered files that feature the same particular number of storm cells. In \cite{HentzenThesis} we describe this parameter study in detail.

%\paragraph{Dynamic clustering}
%Storm cells that are close to each other and move in similar directions for the next ten minutes are grouped into clusters. 10 minutes later, when the new nowcast with updated information about the storm cells comes in, the clustering is updated as well. Instead of keeping the clustering ststic over the horizon of the simulation, this \textit{dynamic} clustering approach allows to account for the merging and splitting of storm cells, as well as the appearance of new storm cells.

%Figure \ref{fig:cluster_ex} shows an example for the application of the dynamic k-mean clustering to the observations present in the available nowcast data at time 00:00 UTC and 30 minutes later at time 00:30 UTC. While clusters 1,5 and 8 stay more or less static over the analyzed 30 minute period, some elements of cluster 3 merge with cluster 2 as the storms move closer to each other. Cluster 4 splits up in two parts, creating 4' and 6'. Clusters 6 and 7 approach each other and merge into 7'. During the 30 minute period, two new storms have appeared and are clustered into one of the eight clusters. Note that in a scenario where the difference in storm cells from one timestep to the next is significant, the number of clusters also adapts in accordance with the look-up table described above.

%\begin{figure}
%	\centering
%	\includegraphics[width= 0.9\columnwidth]{figures/cluster_ex.pdf}
%	
%	\caption{K-mean clustering (K=8) applied to storm cells present in nowcast on December 19th 2016, 00:00 UTC and 00:30 UTC.}
%	\label{fig:cluster_ex}
%\end{figure}

Let $\Xi_{\tau} = \{{C^k_{\tau}}\}$ for $k=1,\dots,K$ and $\tau = 1,\dots,N$ be the new obstacle set containing the K elliptical \textst{obstacles} unsafe-to-fly regions forecast for time $t_0 + \tau$ by the nowcast published at time $t_0$. We now require a stochastic representation of the obstacle set, expressed by $p_{\Xi_\tau}(x,y)$. Recall that this function returns the probability that the coordinate $(x,y)$ is within any obstacle cluster $C^k_{\tau}$ ($k \in {1,2,\dots,K}$)  in $\Xi_\tau$ at time $t_0 + \tau$. 

\subsection{Stochastic Representation of Unsafe-to-fly Region Clusters}
\label{subsec:CovFunc}

\begin{figure}[t!]
	\centering
	\includegraphics[width= 1\columnwidth]{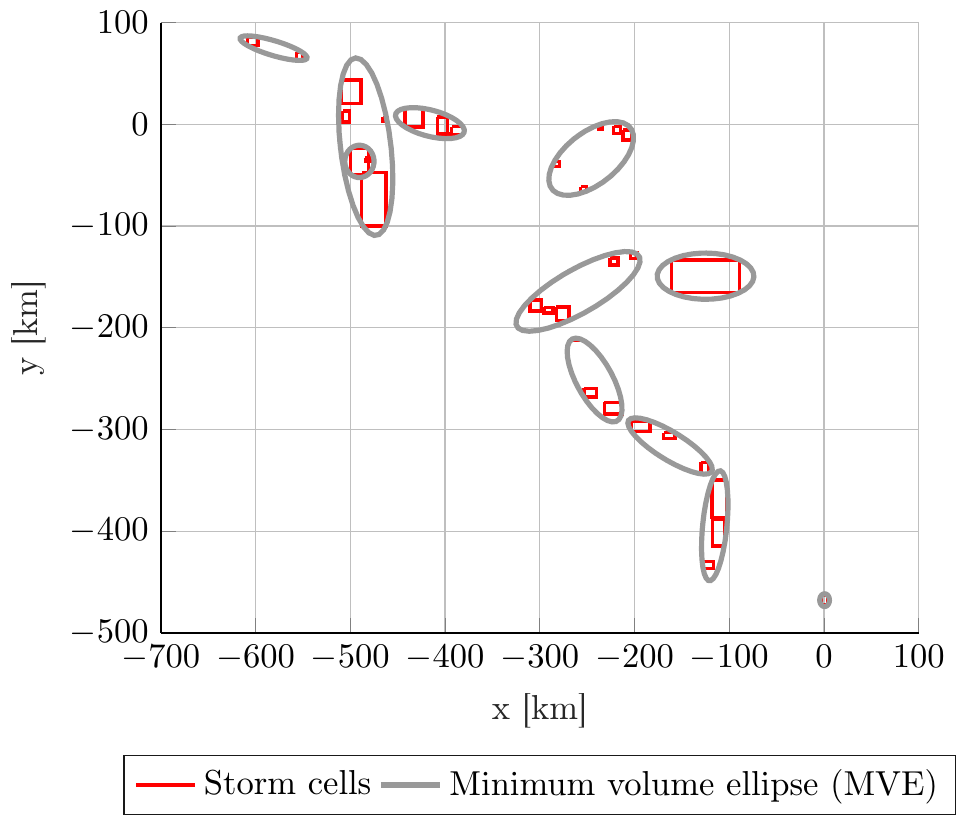}
	
	\caption{Minimum volume ellipse enclosing of clusters in nowcast data of December 19th 2016, 10:30 UTC}
	\label{fig:cluster_ex}
\end{figure}

The computation of $p_{\Xi_\tau}(x,y)$ requires a prediction of the obstacle clusters' position at forecast horizons $\tau$. This prediction in turn calls for a geometric description of the clusters. For this purpose, we choose minimum-volume ellipses (MVE) that enclose the extremities of the thunderstorm cells present within one cluster $C^k_{\tau}$. Ellipses are preferred over polygons for the geometric obstacle characterization due to their efficient mathematical description. Figure \ref{fig:cluster_ex} shows an example of MVE enclosing for two clusters.

Each ellipse $\varepsilon_\tau^k$ is  parameterized by its center $\mathbf{m}_\tau^k \in \R^2$ and its positive definite eccentricity matrix $M_\tau^k$. The position given by coordinate $\mathbf{x}=(x,y)$  is within an elliptic \textst{obstacle} unsafe-to-fly region at time $t_0 + \tau$ if the following condition is fulfilled for any $k$:
\begin{equation*}
\mathbf{x} \in \varepsilon_\tau^k(\mathbf{m}_\tau^k,M_\tau^k) \iff (\mathbf{x}-\mathbf{m}_\tau^k)^T M_\tau^k (\mathbf{x}-\mathbf{m}_\tau^k) \leq 1 .
\end{equation*}
The objective is to compute $p_{\Xi_\tau}$ for the defined obstacle set at forecast horizons $\tau=1,\dots,N$, that is, the probability
\begin{equation*}
p_{\Xi_\tau}=P\{\mathbf{x} \in \varepsilon_\tau^k(m_\tau^k,M_\tau^k) \} \quad \text{for some} \quad k \in \{1,\dots,K\}.
\end{equation*}
%\textcolor{red}{In the example given in \cite{Summers2011}, the error on the ellipse center $m^k$ and the rotation angle of the ellipse $\alpha$ were the parameters used for the stochastic characterization of the random obstacle set. However, this representation fails to account for the change in width and height of the individual storm cells and thus for the change in half-axis length of the resulting bounding ellipse. Moreover, the relation between the individual forecast errors of the storm cells within the ellipse and the resulting deviation between the ellipse bounding the forecast storm cells and the ellipse bounding the true observations of the storm cells is hard to express as a mathematical function and depends on the exact configuration of the storm cells within the ellipse. Thus, the error on the ellipse center $m_k$ cannot be used as a stochastic parameter.}
%One option would be to express the ellipse parameters $m_k^k$ and $M_k^k$ in a stochastic way, as we have done for a simple, single-storm example in \cite{Hentzen2017}. By sampling different realizations of the ellipse parameters, the desired probability function $p_{\Xi_t}$ could then be derived using Monte Carlo methods. In a scenario with multiple individual storm cells within one ellipse, it is hard to find a relation between the stochastic dynamics of the storm cells within the ellipse and the stochastic dynamics of the ellipse's parameters. 

%\subsection{MVE Sampling algorithm}
\label{subsec:MVE_sampling}
\begin{figure}[t!]
	\centering
	\includegraphics[width= 1\columnwidth]{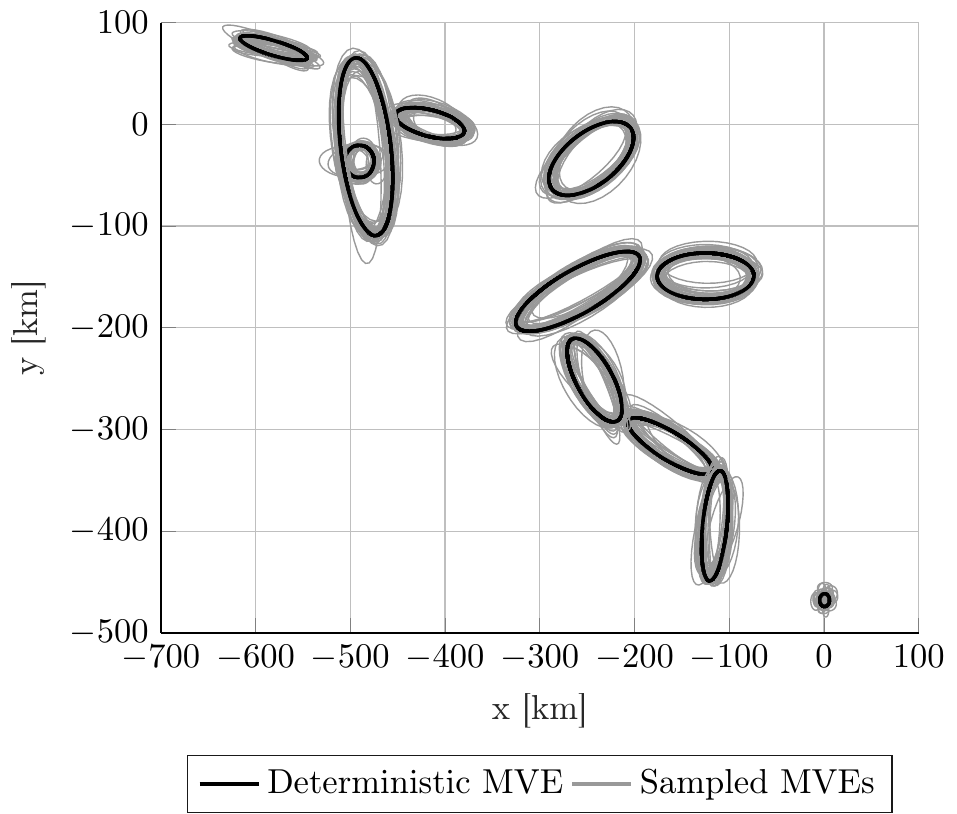}
	
	\caption{50 MVE samples for forecast horizon 10 min.}
	\label{fig:MVE_sampling}
\end{figure}
To achieve this objective we introduce a sampling-based approach, here referred to as the \textit{Minimum-Volume-Ellipse (MVE)} sampling. Every time a new nowcast file is published (at time instance $t_0$), clustering is applied for $t=t_0$ (the current observations of the thunderstorms) and for $t=t_0+\tau$ with $\tau = 1, \dots, N$ (the thunderstorm's forecast at forecast horizons $\tau$). Then, for each of the $K$ clusters present at time $t + \tau$, we predict the position and size of all thunderstorms assigned to that cluster at forecast horizon $\tau$ using the stochastic thunderstorm cell model in subsection \ref{para:stoch_model}. We then enclose the predicted thunderstorms with K MVEs. This procedure is repeated $n_s$ times at each forecast horizon $\tau$, resulting in $n_s$ different MVEs per obstacle cluster. We can then apply Monte Carlo methods to compute the probability $p_{C_\tau^k}(x,y)$ of coordinate $(x,y)\in  X \times Y$ being within the obstacle cluster $C_\tau^k$ at time $t_0 + \tau$.

Figure \ref{fig:MVE_sampling} shows an example for such MVE realizations for a forecast horizon index $\tau$ equivalent to 10 minutes. The black ellipses bound the clustered thunderstorms as provided by the deterministic (uncorrected) data in the nowcast. For each cluster, $n_s = 50$ realizations of the thunderstorm development up to the forecast horizon of 10 minutes are sampled and bounded by an MVE (depicted in grey). This allows to account for the uncertain evolution of the thunderstorms over the 10 minute horizon.

Once the \textit{MVE Sampling} algorithm has been applied to all $K$ obstacle clusters within a nowcast for all forecast horizons $\tau=1,2,\dots,N$, the $K$ individual functions $p_{C^k_\tau}$ for $k=1,\dots,K$ at forecast horizon $\tau$ need to be  merged into the function $p_{\Xi_\tau}$ that returns the probability of $(x,y)$ being in any of the $K$ clusters. We can now approximate the  probability of a coordinate $(x,y)$ being in any of the obstacle clusters as
\begin{equation}
\label{eq:prob_event}
p_{\Xi_\tau}(x,y) \approx \underbrace{1-\underbrace{\prod_{k=1}^K (1-p_{C^k_\tau}(x,y))}_{\text{probability of no thunderstorm at $(x,y)$}}}_{\text{probability of at least one thunderstorm at $(x,y)$}}.
\end{equation}
 In the approximation above, we assume that the event of a coordinate $(x,y)$ falling into a certain cluster is independent of the event of the same coordinate falling into another cluster. While the event that a storm cell is present in a certain region and travels in a certain direction is certainly not independent from the presence of other storm cells in the same region and their movements, we minimize this dependence within the clusters using K-means. Indeed, K-means maximizes the between-cluster variance of all the storm observation vectors, that is, $\mathbf{o} =(x^c, y^c, \xi_{\tau})$. Furthermore, the above inequality is safely conservative. In particular,  denoting stochastic events of $(x,y)$ falling into clusters $C^k_\tau$ and $C^l_\tau$, for $k, l \in \{ 1, \dots, K\}$, with $A$ and $B$ respectively, we have $P(A \cap B) = P(A | B) P(B)$. Since in this setting it is more plausible that $P(A|B) \geq P(A)$ (knowing a cell belongs to one storm cluster increases its probability of belonging to the other storm cluster at the same time), we conclude that $P(A \cap B) \geq P(A) P(B)$. Hence, $1-P(A \cap B) \leq 1-P(A) P(B)$ and $p_{\Xi_\tau} \geq p_{\Xi_{\tau, \text{real}}}$, where $p_{\Xi_{\tau, \text{real}}}$ is the true (unknown) probability. Therefore, our approximation yields an upper bound on the true probability of coordinate $(x,y)$ being in any of the obstacle clusters. This will lead to a more conservative solutions but will not degrade safety.

\begin{figure}[t!]
	\centering
		\includegraphics[width=1\linewidth]{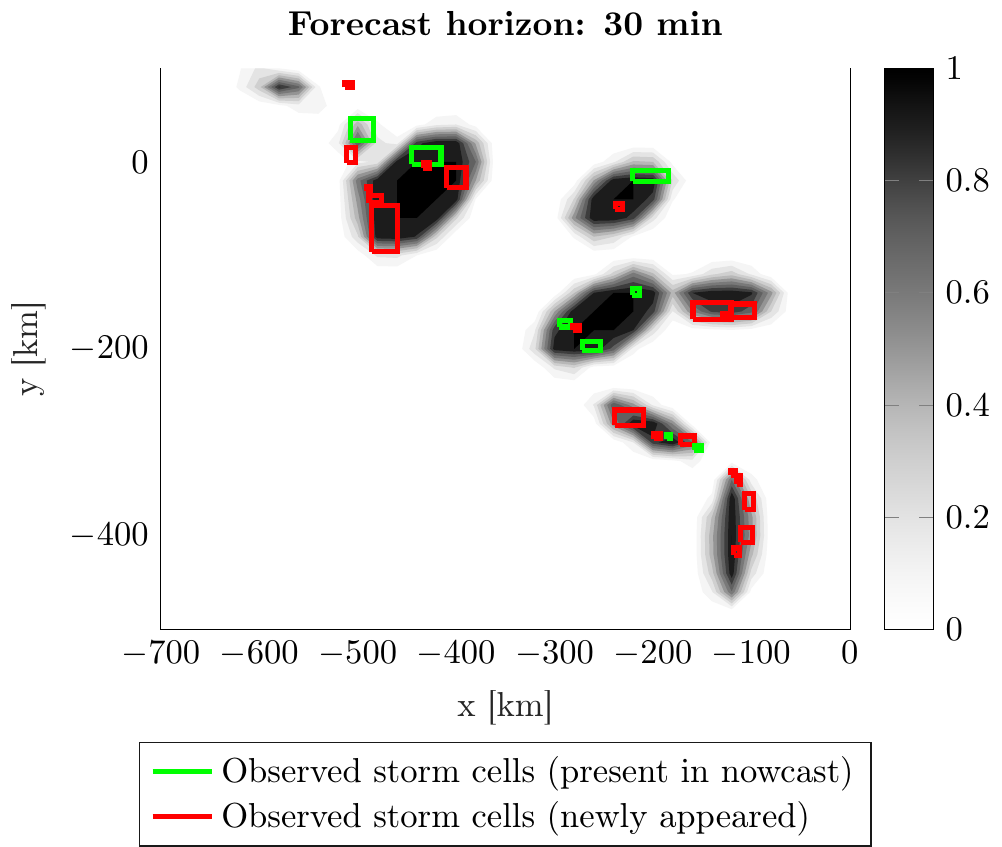}
\caption{Probability of storm $p_{\Xi_{30}}$ for nowcast published on Dec. 19th 2016 at 00:00 UTC.}
	\label{fig:covfct20}
\end{figure}

Figure \ref{fig:covfct20} shows an example of the probability function $p_{\Xi_\tau}(x,y)$ that results for the nowcast provided at 00:00 UTC on December 19th 2016, at a forecast horizon of 30 minutes. The pink, dark boxes indicate thunderstorm cells that were predicted in the original nowcast and indeed observed 30 minutes after the nowcast was issued. The green, light boxes indicate thunderstorm cells that were observed 30 minutes after the nowcast was issued, but not predicted in the nowcast. The figure illustrates the predictive performance of the developed algorithm and the robustness of the clustering approach to newly appearing thunderstorms that are not predicted in the nowcast.

%% file: 3_stoch_reach.tex
% !TeX spellcheck = en_US
% !TeX encoding = UTF-8
% !TeX root = 0_main.tex

%%%%%%%%%%%%%%%%%%%%%%%%%%%%%%%%%%%%%%%%%%%%%%%%
\section{Stochastic Reach-Avoid for Aircraft Trajectory Planning}
\label{sect:stoch_reach}
%%%%%%%%%%%%%%%%%%%%%%%%%%%%%%%%%%%%%%%%%%%%%%%%

\subsection{Stochastic Dynamic Model}
\label{subsec:AircraftDyn}
For the purpose of trajectory planning, we use a simplified model of aircraft dynamics, namely a simple point-mass unicycle model. Furthermore, we assume that the aircraft is flying at constant altitude and constant airspeed $v$. 
Let $S=X \times Y \times [-\pi,\pi] \subset \R^3$ be the state space of the system with state $\mathbf{s}= [x,y,\lambda]\in S$, where $(x, y) \in X \times Y \subset \R^2$ is the two-dimensional position of the aircraft and $\lambda \in [-\pi,\pi]$ is its heading angle measured counter-clockwise from the geographic East. Consistent with current air traffic procedures, we consider  three possible advisories: left turn, right turn and straight flight. Hence, we can represent the control space as  $U = \{-\Omega,0,\Omega\}$, with control input $u \in U$ representing the yaw rate. The dynamics of the aircraft in discrete time are given by
\begin{align}
\label{eq:aircraft_dyn}
\begin{split}
&x_{t+1} = x_t + \delta \cdot (v \cos(\lambda_t) + w_{u,t} )+ \omega^1_t,\\
&y_{t+1} = y_t + \delta \cdot (v \sin(\lambda_t) + w_{v,t} ) + \omega^2_t,\\
&\lambda_{t+1} = \lambda_t + \delta \cdot u_t  + \omega_t^3,\\
\end{split}
\end{align}
where $\delta$ is the sampling time used to discretize the continuous-time model, $w_{u,t}$ and $w_{v,t}$ denote zonal (West-to-East) and meridional (South-to-North) wind values, respectively, and $\mathbf{\omega}_t=[\omega^1_t, \omega^2_t, \omega^3_t]$ is the  disturbance of the system due to wind and model uncertainties. The noise component $\omega^3$  accounts for possible actuation inaccuracies and the fact that the heading cannot be adjusted by the pilots in arbitrary small increments.

%\subsection{Stochastic Reachability}
Given the stochasticity introduced by the disturbance term, the aircraft state at time $t+1$ can be equivalently captured using a so-called stochastic kernel $\chi$. Let $\m{B}(S)$ denote the set of all measurable subsets of the state space $S$. The stochastic kernel $\chi$ returns the probability of state $\mathbf{s}_{t+1}$ being in subset $B \in \mathcal{B}(S)$, given the current state $\mathbf{s}_t$ and the control input $u_t$. The discrete-time aircraft dynamics can thus equivalently be described by
\begin{align*}
\mathbf{s}_{t+1} \sim \chi(B | \mathbf{s}_t, u_t), \quad t = 0, 1, \dots ,N.
\end{align*}
We can state in familiar terms that the state $\mathbf{s}_{t+1}$ is sampled from the probability distribution $\chi(\cdot | \mathbf{s}_t, u_t)$. The numerical computation of $\chi$ will be explained later in this section.

The unsafe-to-fly regions to be avoided by the aircraft are hazardous thunderstorm regions. A stochastic thunderstorm  model was derived in Section \ref{sect:stoch_mod}, which provided  $p_{\Xi_\tau}(x,y)=P\{(x,y) \in \Xi_\tau\}$, the probability of coordinate $(x,y) \in X \times Y$ belonging to a thunderstorm cell at a given time $t_0 + \tau$. The probability of the state $\mathbf{s}=[x,y, \lambda]$ intersecting a thunderstorm can thus be computed as $p_{\Xi_\tau}(\mathbf{s})=p_{\Xi_\tau}(x,y)$ for each $\lambda \in [-\pi, \pi]$.

\subsection{The Reach-Avoid Problem and Solution}
\label{subsec:ProbForm}

Let the sequence $\mu = \{\mu_0,\mu_1,\dots,\mu_{N-1}\}$ be a control policy, that is, a map $\mu_t: S \rightarrow U$ that returns a control input $u_t\in U$ for every state $\mathbf{s}_t \in S$ for $t=0,\dots,N-1$. Given the stochastic kernel $\chi$ and the thunderstorm probability function $p_{\Xi_\tau}$ the reach-avoid objective is as follows. Find  an optimal control policy that maximizes the probability of the aircraft reaching a goal set $G \subset S$ in $N$ steps, while staying in a safe set $Z_t$ for all time steps $t=0, 1, \dots, N$. While the initial obstacle set $\Xi_0$ is assumed to be known and deterministic, the future obstacles are predicted using the parameterized stochastic thunderstorm cell model introduced in Section \ref{sect:stoch_mod}.

The reach-avoid probability for a given control policy $\mu$ and initial state $\mathbf{s}_0, \Xi_0$ is written mathematically as\footnote{The notation implies requirement on $i$ is fulfilled when $\mathbf{s_0} \in G$.}
\begin{align*}
r_{\mathbf{s}_0,\Xi_0}(\mu) \coloneqq P_{\mathbf{s}_0,\Xi_0}^\mu \{&\exists j \in [0,N]: \mathbf{s}_j \in G \\ &\wedge \forall i \in [0,j-1] : \mathbf{s}_i \in Z_i \setminus G \}.
\end{align*}
Let $\mathbf{1}_A$ denote the indicator function of the set $A$. It can be verified that 
\begin{align*}
&r_{\mathbf{s}_0,\Xi_0}(\mu)=\sum_{j=0}^N\bigg(\prod_{i=0}^{j-1} \mathbf{1}_{Z_i\setminus G}(\mathbf{s}_i)\bigg) \mathbf{1}_G(\mathbf{s}_j)
 \\
&=\begin{cases*}
1 & if $\exists j \in [0,N] : \mathbf{s}_j \in G \wedge \forall i \in [0,j-1] \: \mathbf{s}_i \in Z_i \setminus G$ \\
0        & otherwise.
\end{cases*}
\end{align*}

Denoting by $p_{Z_i \setminus G}(\mathbf{s}_i) = \mathbf{1}_{Z\setminus G}(\mathbf{s}_i) - p_{\Xi_i}(\mathbf{s}_i)$ the probability of being in the safe set but not in the goal set, the above can equivalently be written as (see \cite{Summers2011} for a proof)  
\begin{align*}
\begin{split}
r_{\mathbf{s}_0,\Xi_0}(\mu)&= E_{\mathbf{s}_0}^\mu \bigg[\sum_{j=0}^N\bigg(\prod_{i=0}^{j-1} p_{Z_i \setminus G}(\mathbf{s}_i)\bigg)\mathbf{1}_G(\mathbf{s}_j)\bigg].\\
%& := V_0^\mu(x_0)
\end{split}
\end{align*}

To maximize the probability of safely reaching the target, we thus need to compute $r^*(\mathbf{s}_0,\Xi_0) : = \max_{\mu} r_{\mathbf{s}_0,\Xi_0}(\mu)$. We achieve this by solving the stochastic optimal control problem via dynamic programming. 

We define the optimal cost-to-go function $V$$_t^* : S\rightarrow [0,1], t=0,\dots,N$. It maximizes $V_0^\mu(\mathbf{s}_0) := r_{\mathbf{s}_0,\Xi_0}(\mu)$. The backwards recursion needed to compute the optimal cost-to-go $V_t^{*}$ for $t=N-1,\dots,0$ is based on the dynamic programming principle of optimality. In particular,
\begin{align}
\label{eq:OptimalCost}
V_{N}^{*}(\mathbf{s}) = &\mathbf{1}_G(\mathbf{s})\\
\nonumber
V_{t}^{*}(\mathbf{s}) = &\mathbf{1}_G(\mathbf{s}) 
+ p_{Z_t\setminus G}(\mathbf{s})\max_{u_t \in U} \int_{S} V^{*}_{t+1}(\mathbf{s}') \cdot \chi(\mathbf{s}'| \mathbf{s},u).
\end{align}
The corresponding optimal reach-avoid policy $\mu^*_t$ for $t=0,\dots,N-1$ is given by
\begin{align*}
%\begin{split}
%\label{eq:OptimalControl}
\mu_t^*(\mathbf{s}) \in &\arg \max_{u \in U}\bigg\{\mathbf{1}_G(\mathbf{s}) + p_{Z_t\setminus G}(\mathbf{s}) \\
&\int_{S} V^{*}_{t+1}(\mathbf{s}') \cdot \chi(\mathbf{s}'| \mathbf{s}, u)\bigg\}.
%\end{split}
\end{align*}

\subsection{Numerical Implementation}
\label{subsec:num_impl}

To compute the above reachability probability, first we discretize the state space and compute the transition  probability matrix $p_r(\mathbf{s}'| \mathbf{s},u) : S \rightarrow [0,1]$, where $\mathbf{s},\mathbf{s}' \in S, u \in U$ denote state and input values at discretized points in the respective domains. The function $p_r(\mathbf{s}'| \mathbf{s},u)$ is derived based on the transition kernel $\chi$ and maps to each (discretized) state $\mathbf{s}'$ the probability of being reached at the next time step, given the current discretized state $\mathbf{s}$ and control input values $u$. Next, we compute recursion \eqref{eq:OptimalCost} at discretized grid points using matrix multiplication and array maximization. 

The state space $S$ is discretized in all three dimensions $X$, $Y$, and $[-\pi,\pi]$ using the discretization sizes $dx$, $dy$ and $d\lambda$, respectively. The process disturbance $\omega$ is drawn from a normal distribution $\m{N}$ with zero mean and covariance matrix $\Sigma_\omega$ (\cite{Summers2011}):
\begin{equation*}
\mathbf{\omega}=[\omega^1, \omega^2, \omega^3]^T \sim \m{N}(0,\Sigma_\omega),
\end{equation*}
where $\Sigma_{\mathbf{\omega}}= \text{diag}(\sigma_1^2,\sigma_2^2,\sigma_3^2)$, ($\sigma$ denoting a standard deviation). Given the zero mean uncertainty, the mean state $\mathbf{s}'_m$ that the aircraft will reach at the next time step given its current state $\mathbf{s} = [x,y,\lambda]$ and the control input $u$ is 
\begin{equation*}
\mathbf{s}'_m =
\begin{bmatrix}
x + \delta \cdot v \cos(\lambda) + \delta \cdot w_{u}\\
y + \delta \cdot v \sin(\lambda) + \delta \cdot w_v\\
\lambda  + \delta \cdot u\\
\end{bmatrix}.
\end{equation*}
Using the probability density function of the  multivariate normal distribution $f_n$ with mean $\mathbf{s}'_{m}$ and covariance $\Sigma_{\mathbf{\omega}}$,  
we obtain the discretized transition probability $p_r$ in two steps:
\begin{equation*}
p_{r,temp}(\mathbf{s}' | \mathbf{s},u) = dx \cdot dy \cdot d\lambda \cdot f_n(\mathbf{s}'_m ),
\end{equation*}
\begin{equation*}
p_{r}(\mathbf{s}' | \mathbf{s},u) = \frac{p_{r,temp}(\mathbf{s}' | \mathbf{s},u)}{\sum_{S}p_{r,temp}(\mathbf{s}' | \mathbf{s},u)}.
\end{equation*}
This is then used to compute the integral term in \eqref{eq:OptimalCost} as a summation over the state 
\begin{align*}
\int_{S} V^{*}_{t+1}(\mathbf{s}') \cdot \chi(\mathbf{s}'| \mathbf{s}, u) \approx \sum_S V^{*}_{t+1}(\mathbf{s}') \cdot p_r(\mathbf{s}'| \mathbf{s}, u). 
\end{align*}
The other terms in computation of \eqref{eq:OptimalCost} are simply evaluations of a point belonging to a set or maximization over an array of numbers. Note that finer discretization grids result in larger matrix and array dimensions and hence more computational effort. 

%% file: 4_case_studies.tex
% !TeX spellcheck = en_US
% !TeX encoding = UTF-8
% !TeX root = 0_main.tex

%%%%%%%%%%%%%%%%%%%%%%%%%%%%%%%%%%%%%%%%%%%%%%%%
\section{Case Studies}
\label{sect:case_studies}
%%%%%%%%%%%%%%%%%%%%%%%%%%%%%%%%%%%%%%%%%%%%%%%%

The developed trajectory planning tool is applied to two realistic flight scenarios using state-of-the-art nowcast data. The safety of the trajectories optimized based on stochastic control framework are compared against the actual storm observations. 

The proposed algorithm is implemented in Matlab and simulated on an Intel Core i7 2.5 GHz processor with 16 GB RAM. All CPU times refer to this infrastructure.

\subsection{Problem Setup}

\begin{figure}[t!]
	\centering
	\includegraphics[width=0.9\linewidth]{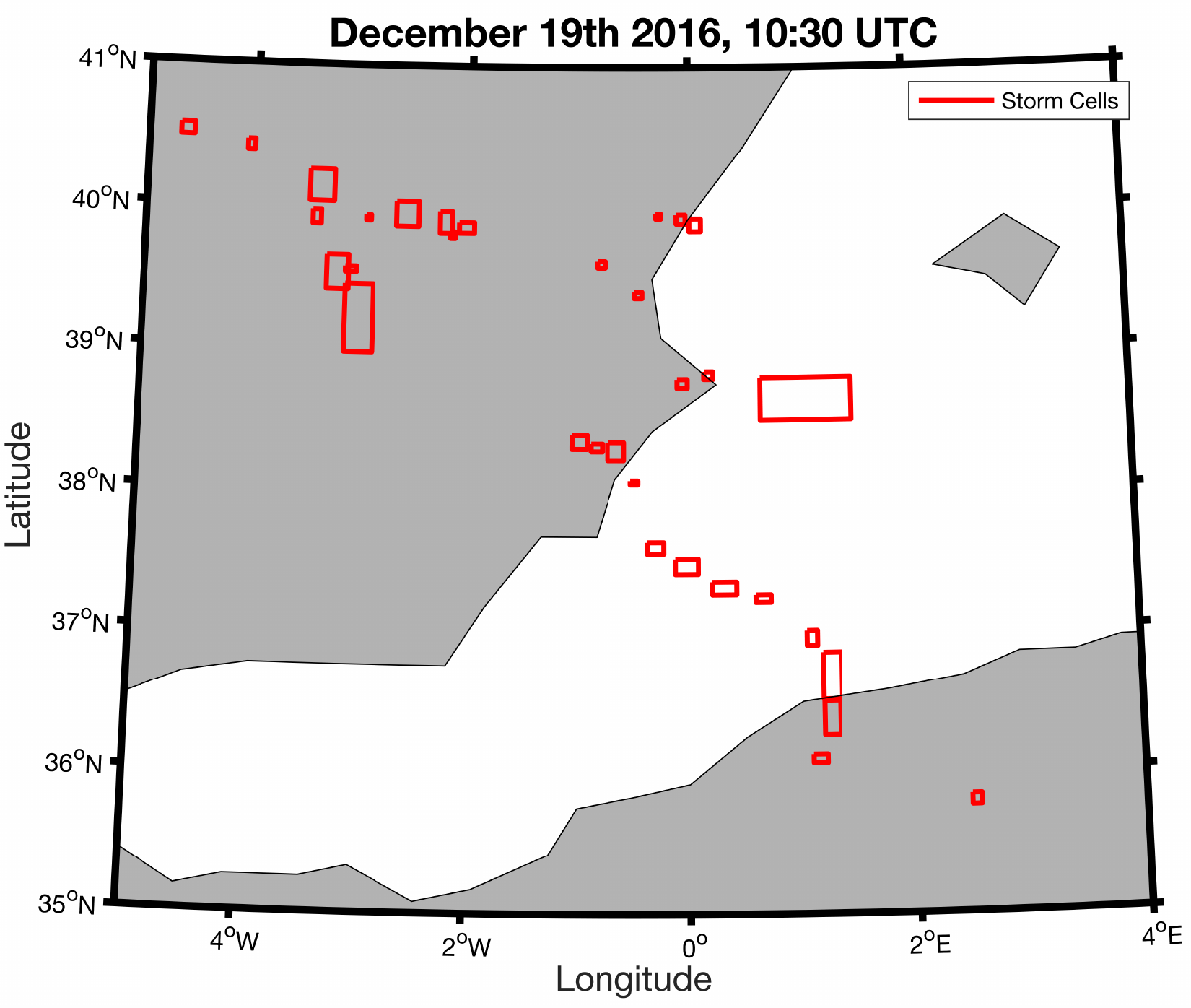}
	\caption{Illustration of storm cells in nowcast on Dec. 19th 2016, 00:00 UTC.}
	\label{fig:nowcastcs}
\end{figure}

%\paragraph{Flight Details}
A region of airspace over the Spanish eastern coast and the Balearic sea, spanned by the corner coordinates $(-5\degree,41\degree)$ and $(4\degree,35\degree)$ is considered. This region is projected onto a 2D plane with a Lambert conformal projection centered at $(38^\circ, -98^\circ)$. The area of interest in $X\times Y$ coordinate is $[-650,100]\times[-550,200]$ (in km).  To simulate a realistic scenario, the storm cell configuration is extracted from the AEMET nowcast of December 19th 2016 at 10:30 UTC, illustrated in Figure \ref{fig:nowcastcs}.

\emph{Flight 1:}
The initial condition of the aircraft is  $\textbf{s}_0= (-365,100, -1.92$ \text{rad}) heading directly to its goal set $G = [-490,-470]\times[-230,-210]\times[-\pi,\pi]$. From this starting position, the goal set is at 340 km of the starting point (26 min straight line flight at the chosen airspeed, without considering wind).

\emph{Flight 2:}
The aircraft starts at  $\textbf{s}_0 = (-100,30,-2.5$ \text{rad}) heading directly to the same goal set as in Flight 1. In contrast to Flight 1, the aircraft starts further away from the goal set, at a distance of 455 km (35 min straight line flight at the given airspeed, without considering wind).

%\paragraph{Simulation Settings}

We model the aircraft using the discrete time unicycle dynamics presented in the previous Section, with an airspeed of 792 km/h and a maximal turn rate of 0.3 rad/min. As we are dealing with turn advisories, a low turn rate is chosen to  enable a smooth passenger friendly trajectory.  The disturbance parameters in aircraft dynamics in \eqref{eq:aircraft_dyn} are set to  $\Sigma_{\mathbf{\omega}}= \text{diag}(\sigma_1^2,\sigma_2^2,\sigma_3^2)$, with $\sigma_1^2 = \sigma _2^2 = 0.25$ [$m^2$], $\sigma_3 = 4 \cdot 10^{-5}$ [$rad^2$]. Wind forecast values are $w_{u,t}=2.6$ m/s and $w_{v,t}=5.6$ m/s for the zonal (West-to-East direction) and meridional (South-to-North direction) component, respectively. They are assumed to be constant over time and have been retrieved from the corresponding ECMWF Forecast as the average values in the area and time period of interest. At the specified true airspeed $v$ and discretization time step of $\delta = 2$ min, the random variable corresponding to the distance covered by the aircraft in one time step has a 95\% confidence interval $26.4 \pm 1.74$ km. At a yaw rate of $\Omega = 2.7$ rad/min, the random variable corresponding to the turning angle covered in one time step  has a 95\% confidence interval $5.4 \pm 0.014$ rad. 

The time horizon is set to 40 minutes, equivalent to the maximum usable forecast horizon provided in the nowcast files. We note that the chosen time horizon may in theory be extended to an arbitrarily large value to fit any flight duration. However, we are constrained by the maximum forecast horizon provided in the available nowcast files. To illustrate the influence of the time horizon, flight 2 is simulated with both a forecast horizon of 40 minutes and the maximum allowable 60 minutes. Two simulation settings are used: Setting 1 uses a coarse discretization grid and a low number of samples and clusters, while setting 2 uses a very fine grid and a higher number of samples and clusters. Table \ref{tab:parameters_simu} provides an overview of both settings.

\begin{table} [b!]
	%\scriptsize
	\centering
	\begin{tabular}{p{4cm} p{4cm}}
		\toprule
		\textbf{Setting 1 - Coarse} & \textbf{Setting 2 - Fine}\\
		\midrule
		$33 \times 28 \times 32$ grid & $66 \times 56 \times 32$ grid\\
		100 MVE samples & 500 MVE samples\\
		12 clusters & 18 clusters\\
		\bottomrule
	\end{tabular}
	\caption{Overview of the two simulation settings.}
	\label{tab:parameters_simu} 
\end{table}
\subsection{Results}

\begin{figure}[t!]
	\centering
	\includegraphics[width=0.9\columnwidth]{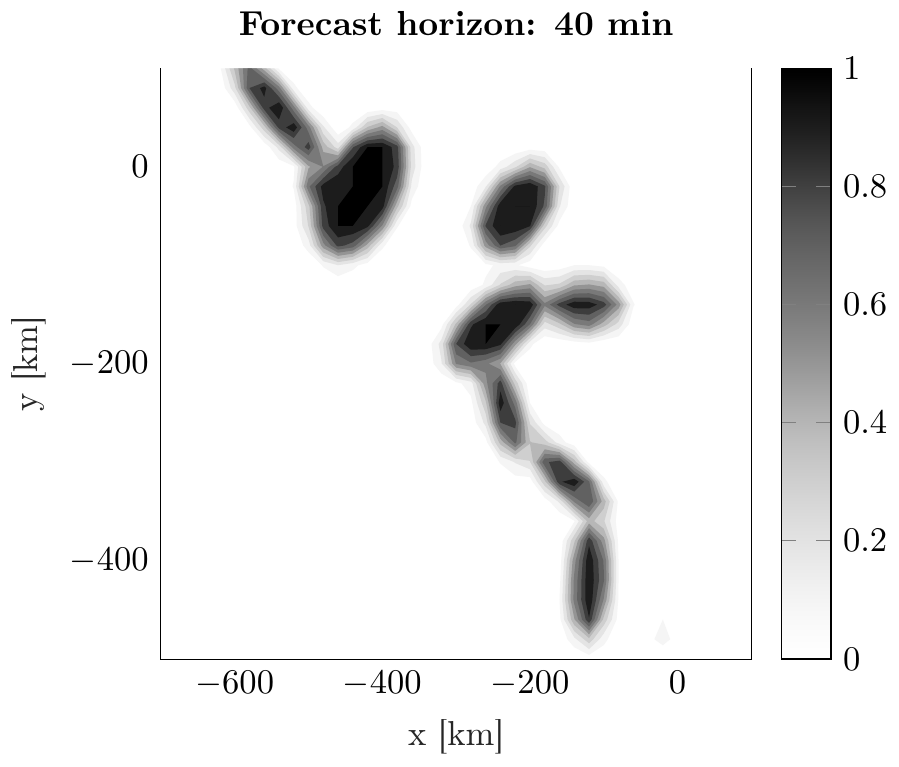}
	\vspace{-0.5cm}
	\caption{Probability of storm $p_{\Xi_{40}}$ computed with the coarse grid setting at forecast horizon of 40 minutes.}
	\label{fig:covfset2_2}
\end{figure}
Figure \ref{fig:covfset2_2} shows the \emph{probability of storm} $p_{\Xi_t}$, in the considered airspace at a 40 minutes forecast horizon, computed with setting 2. The clusters and associated MVEs are clearly distinguishable. While the computation of the storm probability takes 279 s using the coarse setting 1, the simulation with the fine setting 2 significantly degrades computational performance (7'016 s) as seen in Table \ref{tab:comp_perf}.

\begin{table}[b!]
	%\scriptsize
	\centering
	\begin{tabular}{p{5.4 cm} p{1.35cm} p{1.35cm}}
		\toprule
		&\textbf{Setting 1} & \textbf{Setting 2}\\
		\midrule
		Probability of storm $p_{\Xi_\tau}$ ($\tau = 1,\dots,N	$)& 279 s & 7'016 s\\
		Reach-avoid probability $V_0^*(\mathbf{x}_0)$ & 1'196 s & 9'101 s\\
		\bottomrule
	\end{tabular}
	\caption{Overview of computational performance.}
	\label{tab:comp_perf} 
\end{table}

\begin{figure}[t!]
	\centering
	\includegraphics[width=0.9\columnwidth]{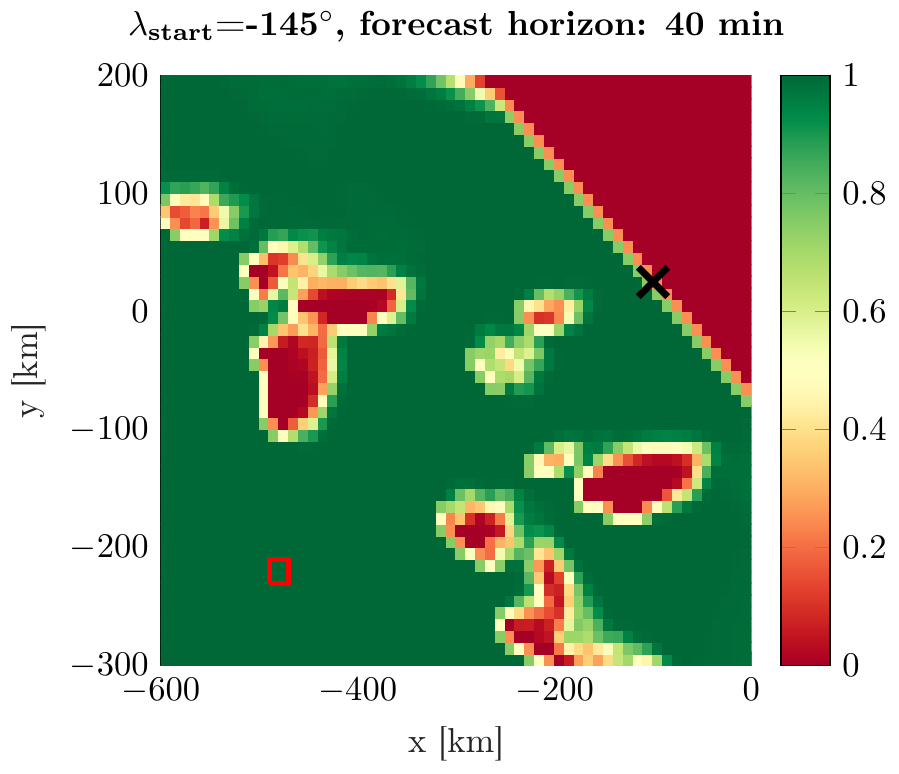}
	\vspace{3mm}
	\includegraphics[width=0.9\columnwidth]{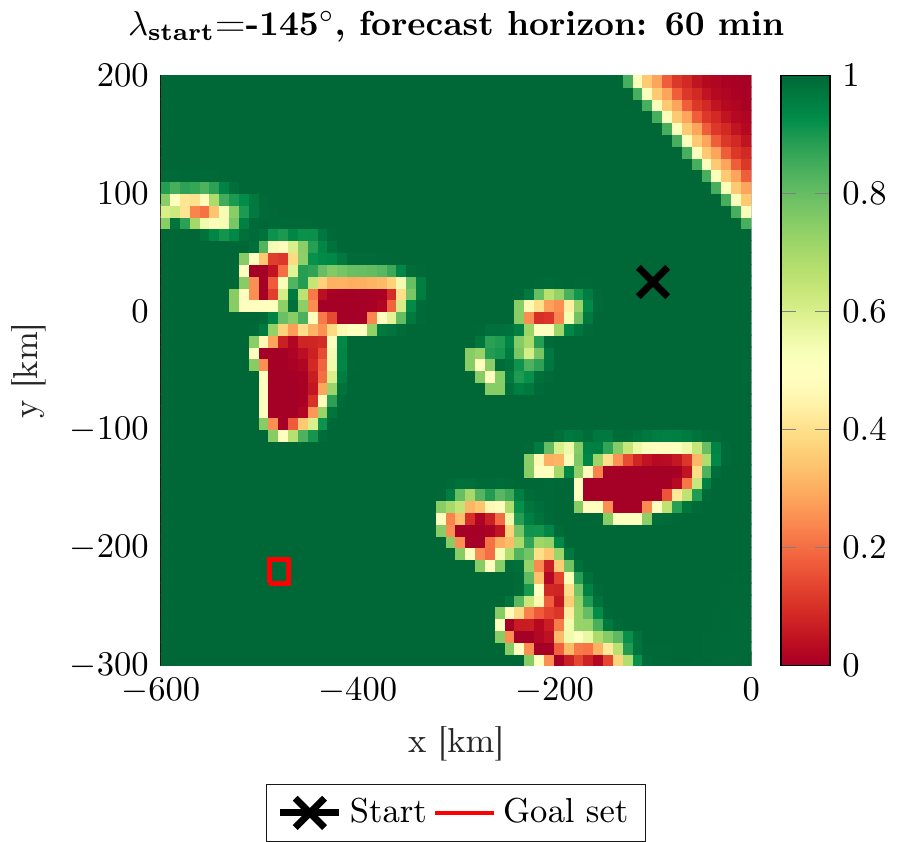}
	\caption{Maximum probability of reaching the desired waypoint while avoiding unsafe-fly regions,  $V_0^*(\mathbf{s}_0)$ computed for flight 2, using two different forecast horizons.}
	\label{fig:valfcts1s2}
\end{figure}

Figure \ref{fig:valfcts1s2} shows the \emph{optimal reach-avoid probability} for flight 2. For every point of the map, this is the maximum probability of reaching the goal set safely. Next to the low probability areas that are due to likely storm presence, we observe the reachability limit in the form of a distinct boundary in the upper right corner. Indeed, given the fixed time horizon and the constant airspeed of the aircraft, the aircraft can only reach targets within a given distance starting with certain heading angle. We recompute the optimal reach-avoid probability for an extended time horizon of 60 minutes. Figure \ref{fig:valfcts1s2} shows that the mentioned reachability limit in the upper right corner has retreated. The aircraft can now safely reach the given goal set from a greater distance.

Runtime for the computation of the optimal reach-avoid probability increases  with grid size, as shown in Table \ref{tab:comp_perf}.

\begin{table*}[t!]
	%\scriptsize
	\centering
	\begin{tabular}{p {3.5cm} p{1.25cm} p{1.25cm} p{1.25cm}p{1.25cm}}
		\toprule
		&\textbf{Flight 1} & & \textbf{Flight 2} & \\
		\midrule
		Setting & 1 & 2 & 2 & 2 \\
		Time horizon (min) & 40 & 40 & 40 & 60 \\
		Flight time & 1'824 s & 1'776 s & 2'280 s & 2'340 s\\
		Reach-avoid probability & 1.00 & 1.00 & 0.56 & 1.00 \\
		Successful trajectories & 100 \% & 100 \% & 88 \% & 95\%\\
		\bottomrule
	\end{tabular}
	\caption{Overview of simulation results. \vspace{2mm}}
	\label{tab:sim_res} 
\end{table*}

\paragraph{Flight 1; Setting 1}

\begin{figure*}[th!]
	\centering
	\includegraphics[width=0.99\linewidth]{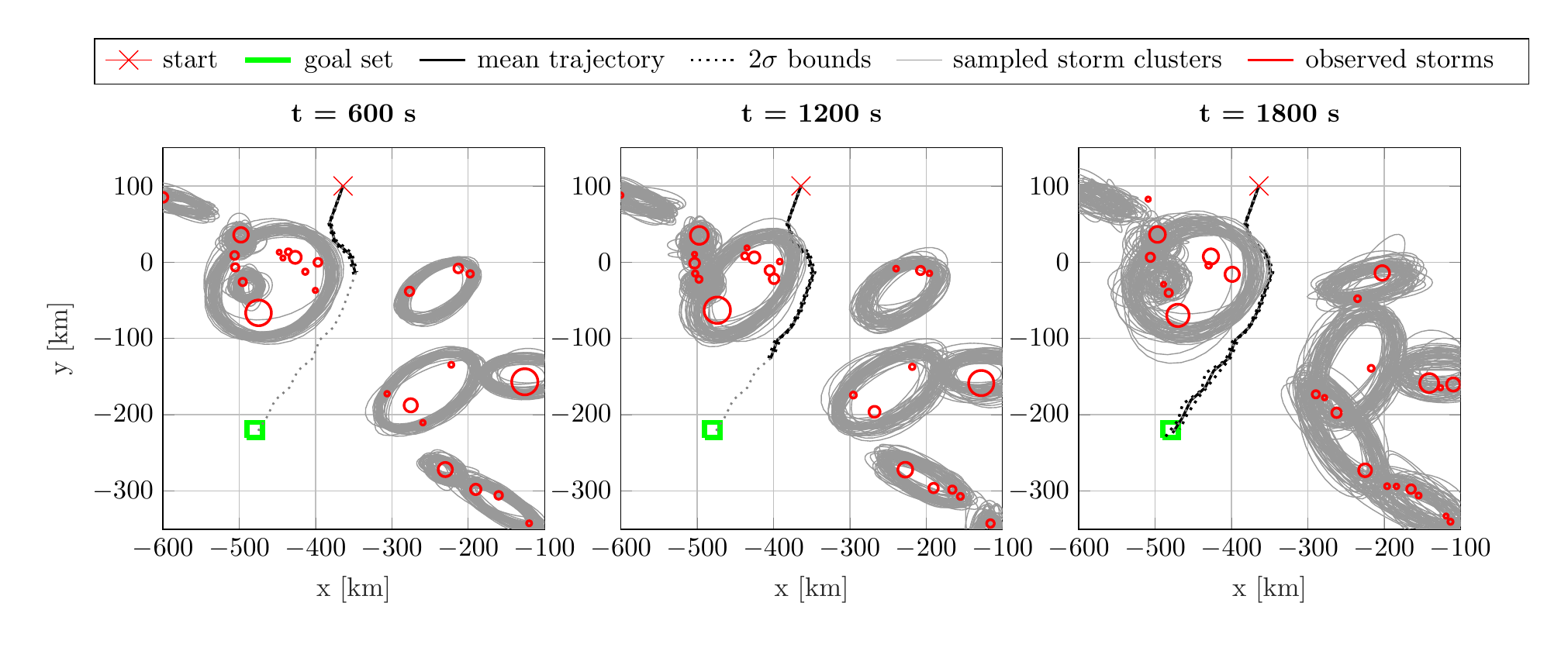}
	\vspace{-0.5cm}
	\caption{Flight 1, Setting 1: mean trajectory and 2$\sigma$ trajectory bounds.}
	\label{fig:trajmve1}
\end{figure*}

Figure \ref{fig:trajmve1} shows the mean trajectory and $2 \sigma$ confidence bounds resulting from 10'000 trajectory simulations for flight 1 with setting 1. The red circles represent the observed storm cells at the given time. The grey ellipses are MVE samples (introduced in Section \ref{subsec:MVE_sampling}), that is, samples of possible future positions of the storm cell clusters. As discussed in Section \ref{subsec:StochNoFly}, the developed algorithm generates trajectories that avoid these areas of high storm probability. This is illustrated in Figure \ref{fig:trajmve1}: Accounting for the inherent uncertainty in the forecasts, the trajectory avoids not only individual predicted storm cells, but entire areas predicted to have high probability of storm. Almost all observed storm cells (depicted in red) fall into one of these areas. In this particular scenario, the trajectory is headed straight towards the goal set at first, following its reference trajectory. To avoid the approaching storm clusters, the aircraft then changes heading and passes through the opening between two clusters. Once the aircraft has passed the clusters, it changes heading once more to head to the goal set. This trajectory, if executed, would avoid all observed storm cells.

The optimal reach-avoid probability for the chosen initial condition is 1.00. All of the simulated trajectories reach the goal set within the specified horizon without hitting one of the observed thunderstorms. Table \ref{tab:sim_res} provides an overview of the relevant simulation results.

\paragraph{Flight 1; Setting 2}
The same flight is simulated with the fine settings. 

Due to the improved resolution (higher number of clusters and smaller discretization size), more regions of the airspace are predicted to be safe. The resulting trajectory is therefore less conservative, passing closer to the storms and reaching the goal set 168 s earlier. The above comparison of a coarse and a fine simulation setting shows that both settings yield a safety-optimal solution. Choosing a finer resolution will potentially yield better performance in terms of flight time (and therefore fuel, $CO_2$ emissions, costs), as the increased accuracy in storm representation reduces the area of unsafe-fly regions. However, computational runtime increases significantly. The choice of resolution thus represents a trade-off between system performance in terms of flight time and computational runtime.

\begin{figure*}[t!]
\centering \includegraphics[width=1\linewidth]{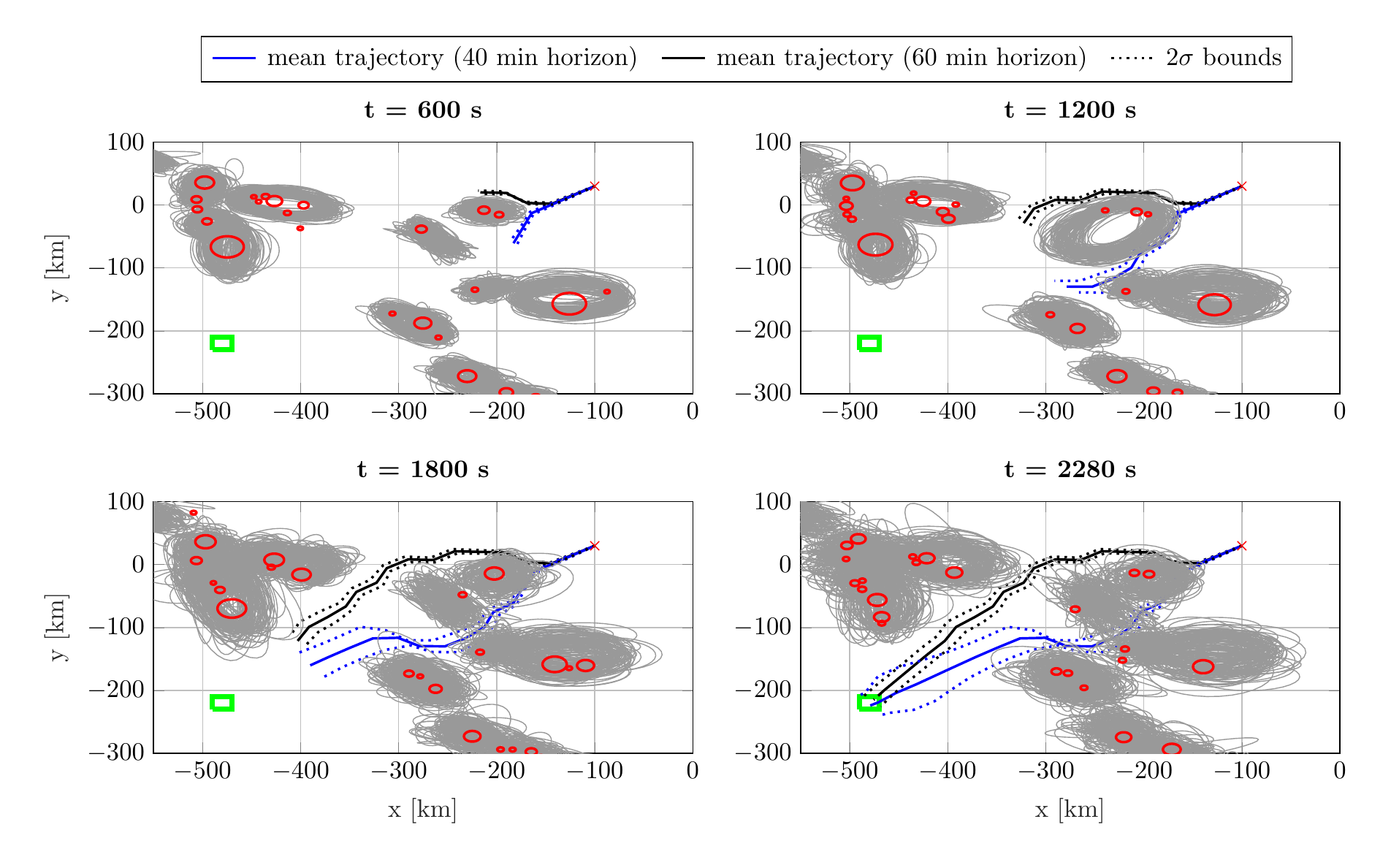}
\vspace{-0.5cm}
\caption{Flight 2, Setting 2: mean trajectory and 2$\sigma$ trajectory bounds.}
\label{fig:trajmve2}
\end{figure*}

\paragraph{Flight 2; Setting 2} For flight 2 we keep the simulation settings 2, but change the initial conditions. Figure \ref{fig:trajmve2} shows the mean trajectory and $2 \sigma$ confidence bounds out of 10'000 simulated realizations for two different time horizons: 40 and 60 minutes. 

Due to the initial conditions on the aircraft state, in particular the long distance to the goal set (considering the limited range of 40 min flight time), the aircraft cannot take a right turn at the start to avoid the dangerous clusters ahead and to minimize the risk of flying through a storm. Instead, the combined reach-avoid probability is maximized by a trajectory that trades storm avoidance for the ability to reach the goal set within the fixed horizon of 40 min. In reality, this would prevent the aircraft from choosing a route that avoids all storms, but fails to reach the target within the given time or fuel constraint, which could result in catastrophic failure. The resulting optimal reach-avoid probability is 0.56. Out of the 10'000 simulated trajectories, 8'773 reach the target. All of them avoid the observed storm cells in red. Some trajectories are forced to pass through forecast storm clusters. However, no storms are eventually observed in these areas. This explains the discrepancy between the value of the reach-avoid probability and the number of successful trajectories. Note that this discrepancy adds conservatism to the predicted safety of the generated trajectory and is not detrimental to safety.

We computed the maximum reach-avoid probability for flight 2 with an extended time horizon of 60 minutes. The aircraft now has the possibility to divert towards north of the predicted storm clusters. While the 40 minute trajectory reaches its goal at 2280 s on average,  the 60 minute trajectory does not get as close to this time constraint, reaching the goal set after 2340 s on average. In contrast to the 40 minute time horizon, in the 60 minute horizon, 94\% of the simulated trajectories reach the goal state while avoiding all observed storm cells. This confirms that the time horizon should be determined according to the available forecast data and underlines the importance of high-quality forecasts with extended forecast times. In practice, however, there may be fuel limitation or time of flight constraints for choosing the maximum time horizon.

%% file: 5_conclusion.tex
% !TeX spellcheck = en_US
% !TeX encoding = UTF-8
% !TeX root = 0_main.tex

%%%%%%%%%%%%%%%%%%%%%%%%%%%%%%%%%%%%%%%%%%%%%%%%
\section{Conclusions}
\label{sect:conclusion}
%%%%%%%%%%%%%%%%%%%%%%%%%%%%%%%%%%%%%%%%%%%%%%%%
We presented a novel approach to extract a  probabilistic dynamic description of hazardous thunderstorm regions from state-of-the-art nowcast data. We integrated the stochastic thunderstorm model in a stochastic, optimal trajectory planning tool, which maximizes the probability of reaching a specified goal set while avoiding the hazardous thunderstorm regions. Simulations of different scenarios using historic nowcast files have proven the obstacle avoidance quality of the generated trajectories, as well as the accuracy of the computed safety levels. This information  could be presented to both pilots and controllers to aid them in decision making processes. 
The computational performance is nevertheless not satisfactory for online implementation. The main reason is the curse of dimensionality inevitably resulting from the use of the dynamic programming approach in a 4D spatiotemporal space. Hence, the developed algorithm  could be used in pre-tactical scenarios and as a benchmark tool to compare trajectories generated by other optimization approaches with the safety-optimal trajectories. To improve runtime, future work should explore alternative, more efficient numerical optimization methods such as nonlinear programming or chance-constrained optimization, and parallelized GPU computing. More efficient computational performance would also enable to explore control frameworks such as receding horizon control. This would enable replanning of the trajectory as new forecasts come in. Another potential field of future research is the incorporation of other adverse meteorological phenomena, e.g., wind-shear, hail, or clean air turbulence, into the presented framework, further maximizing safety.